\documentclass[11pt,reqno]{amsart}

\numberwithin{equation}{section}

\usepackage{amssymb}
\usepackage{amsmath}
\usepackage{amsbsy}
\usepackage{amscd}
\usepackage{amsthm}
\usepackage{amsfonts}
\usepackage{enumerate}

\usepackage{color}

\usepackage{ifpdf}
\ifpdf \usepackage[pdftex,pdfstartview=FitH,pdfpagemode=none,colorlinks,bookmarks,linkcolor=blue]{hyperref} \else  \usepackage[hypertex]{hyperref} \fi

\usepackage[nobysame,abbrev,alphabetic]{amsrefs}
\newcommand{\hide}[1]{\marginpar{Hidden Text}}

\newtheorem{theorem}{Theorem}[section]

\newtheorem{lemma}[theorem]{Lemma}
\newtheorem{corollary}[theorem]{Corollary}
\newtheorem{definition}[theorem]{Definition}

\newtheorem{conjecture}[theorem]{Conjecture}

\newtheorem{proposition}[theorem]{Proposition}

\newtheorem{hypothesis}[theorem]{Hypothesis}

\newtheorem{remark}[theorem]{Remark}
\theoremstyle{definition}

\newcommand{\ksum}{\sum_{k=k_-}^{k_+}}
\newcommand{\kqsum}{\sum_{\substack{k_-\leq k\leq k_+\\ q_{k+1}>q_k^{\frac\tau 2}}}}
\newcommand{\kqprod}{\prod_{\substack{k_-\leq k\leq k_+\\ q_{k+1}>q_k^{\frac\tau 2}}}}

\newcommand{\bn}{\bar{n}}

\newcommand{\bC}{\mathbb{C}}

\newcommand{\bE}{\mathbb{E}}

\newcommand{\bR}{\mathbb{R}}
\newcommand{\bZ}{\mathbb{Z}}

\newcommand{\bN}{\mathbb{N}}

\newcommand{\bT}{\mathbb{T}}

\newcommand{\tphi}{{\tilde{\phi}}}
\newcommand{\tx}{{\tilde{x}}}

\newcommand{\tpsi}{{\tilde{\psi}}}

\newcommand{\hh}{{\hat{h}}}

\newcommand{\rmm}{{\mathrm{m}}}

\newcommand{\bfn}{\mathbf{n}}

\newcommand{\bfr}{\mathbf{r}}

\newcommand{\bfone}{\mathbf{1}}

\newcommand{\Var}{\operatorname{Var}}

\newcommand{\m}{\mathrm{m}}

\newcommand{\Exp}{\mathop{\mathbb{E}}}
\newcommand{\Prob}{\mathop{\mathbb{P}}}

\newcommand{\onto}{\xymatrix{\ar@{>>}[r]&}}
\newcommand{\da}[4]{\xymatrix{#1 \ar@<.5ex>[r]^{#2} \ar@<-.5ex>[r]_{#3} & #4}}

\newcounter{subconst}[subsection]

\newcounter{const}

\newcounter{CONST}

\begin{document}

\title{M\"{o}bius disjointness for non-uniquely ergodic skew products}
\author[Z. Wang]{Zhiren Wang}
\address{ Pennsylvania State University, University Park, PA 16802, USA}
\email{zhirenw@psu.edu}
\setcounter{page}{1}
\begin{abstract}
For $\tau>2$, let $T$ be a $C^\tau$ skew product map of the form $(x+\alpha,y+h(x))$ on $\bT^2$ over a rotation of the circle. We show that if $T$ preserves a measurable section, then it is disjoint to the M\"{o}bius sequence. This in particular implies that any non-uniquely ergodic $C^\tau$ skew product map on $\bT^2$ has a finite index factor that is disjoint to the M\"{o}bius sequence.
\end{abstract}
\maketitle
{\small\tableofcontents}

\section{Introduction}

Let $T$ be a skew product map on $\bT^2$ over a rotation of the circle $\bT^1$. That is, \begin{equation}\label{DynaEq}T(x,y)=(x+\alpha, y+h(x)),\end{equation}
where $h:\bT^1\mapsto\bT^1$ is continuous, and $\alpha\in[0,1)$.

A measurable invariant section of  $T$ is a graph $(x,g(x))$ where $g:\bT^1\mapsto\bT^1$ is measurable, such that $T(x,g(x))$ is still in the graph for Lebesgue almost every $x$.

Our main result is the following Mobius disjointness property:

\begin{theorem}\label{MainThm} Suppose $h$ is $C^\tau$ for some real value $\tau>2$ and $T$ preserves a measurable invariant section. Then for all $(x,y)\in\bT^2$, and all continuous functions $f\in C(\bT^2)$,
\begin{equation}\label{MainEq}\frac1N\sum_{n\leq N}\mu(n)f(T^n(x,y))\rightarrow 0\text{ \rm as }N\rightarrow\infty.\end{equation}\end{theorem}

One important feature of the theorem is that it holds for all $\alpha$, without assuming any Diophantine condition.

By a dichotomy of Furstenberg \cite{F61}*{Lemma 2.1}, if a map $T$ of the form \eqref{DynaEq} is not uniquely ergodic, then for some positive integer $\xi$, the solution $g(x+\alpha)-g(x)=\xi h(x)$ has a measurable solution $g:\bT^1\mapsto\bT^1$. Let $\pi_\xi:\bT^2\mapsto\bT^2$ be the $\xi$-to-one projection $\pi(x,y)=\pi(x,\xi y)$. Then the transform $T_\xi(x,y)=(x+\alpha, y+\xi h(x))$ is a topological factor of $T$ through $\pi$, in other words, $\pi_\xi\circ T=T_\xi\circ\pi_\xi$. One can easily check that the graph $(x,g(x))$ is a measurable invariant section for $T_\xi$. Hence Theorem \ref{MainThm} implies:

\begin{corollary}\label{MainCor}Suppose $h$ is $C^\tau$ for some real value $\tau>2$ and $T$ is not uniquely ergodic. Then there exists $\xi\in\bN$, such that the $\xi$-to-one topological factor $T_\xi$ of $T$ via the projection $\pi_\xi$ satisfies $$\frac1N\sum_{n\leq N}\mu(n)f(T_\xi^n(x,y))\rightarrow 0\text{ \rm as }N\rightarrow\infty.$$\end{corollary}

The M\"obius function is the multiplicative function
$$\mu(n)=\left\{\begin{array}{ll}(-1)^{\text{\# of prime factors of }n},&\text{ if }n\text{ is square free}\\0,&\text{ otherwise.}\end{array}\right.$$

It is expected that the M\"obius function captures much of the randomness in the distribution of prime numbers. This is characertized by Sarnak's Mobius Disjointness Conjecture, which has been the focus of much research in number theory and dynamical systems in recent years.

\begin{conjecture}\label{SarnakConj}(M\"obius Disjointness Conjecture, \cite{S09}) For a continuous transformation $T: X\mapsto X$ with zero topological entropy, where $X$ is a compact metric space and $T$ is continuous, then for all continuous $f\in C(X)$ and all $x\in X$, $\frac1N\sum_{n\leq N}\mu(n)f(T^n(x))\rightarrow 0$.\end{conjecture}

It should be emphasized the conjecture addresses all points $x\in X$, instead of almost every $x$ with respect to a measure.

When $X$ is a single point, Conjecture \ref{SarnakConj} gives the prime number theorem. The prime number theorem in arithmetic progressions corresponds to the case of a rotation on the finite abelian group $\bZ/q\bZ$.  Davenport's classic theorem \cite{D37} that $\frac1N\sum_{n\leq N}\mu(n)e(\alpha n)\rightarrow 0$ is the case of a rotation of the circle.

Researches during the last a few years confirmed numerous cases of Conjecture \ref{SarnakConj}. To list a few: \cites{MR10, G12, GT12, B13a, B13b, BSZ13, KL15, ELD14, MMR14, P15, FKLM15}. Most of these results use Vinogradov's bilinear method, or its  modern variant due to Vaughan  \cite{V77}. The Bourgain-Sarnak-Ziegler criterion, developed in \cite{BSZ13}, offers a different variant that allows to interprete bilinear averages in as ergodic averages in joinings of dynamical systems.

All the special case listed above are regular dynamical systems. A system $(X,T)$ is regular if the ergodic average $\frac1N\sum_{n\leq N}f(T^n(x))$ converges for every $x\in X$ and every continuous function $f$. In particular, all uniquely ergodic systems are regular by Birkhoff's ergodic theorem.

The only irregular dynamical systems for which M\"obius disjointness has been studied come from the family \eqref{DynaEq}. In this family, for generic choices of $\alpha$ that satisfy certain diophantine conditions, a $C^{1+\epsilon}$-differentiable $T$ is still regular, as well as the joinings needed for applying the Bourgain-Sarnak-Ziegler criterion. For such $\alpha$, Conjecture \ref{SarnakConj} was proved by Ku\l{}aga-Przymus and Lema\'nczyk \cite{KL15}.

On the other hand, for Liouville choices of $\alpha$, the map $T$ can be irregular. Such a counterexample was first constructed by Furstenberg in \cite{F61}. Furstenberg's example is real analytic, but it is not hard to modify it to get a counterexample of finite differentiability.

Irregular transformation from the family \ref{DynaEq} were first studied in \cite{LS15} by Liu and Sarnak. They proved Conjecture \ref{SarnakConj} for a class of analytic skew products of the form \eqref{DynaEq} under an additional assumption that $\hh(m)$ decays not too fast. In \cite{W15}, the author proved Conjecture \ref{SarnakConj} for all real analytic maps of the form \eqref{DynaEq}. One main ingredient from \cite{W15} was the use of the estimate of averages of non-pretentious multiplicative functions on typical short intervals by Matom\"aki, Radziw\l\l{} and Tao \cite{MRT15}. The use in \cite{W15} of this new tool was quantitative, and as a consequence analyticity, or at least the weaker condition that $|\hh(m)|\ll e^{-\tau|m|^{\frac12+\epsilon}}$, had to be assumed to control the estimates from \cite{MRT15}.

The main aim of Theorem \ref{MainThm} is to prove Conjecture \ref{SarnakConj} for some irregular dynamical systems of finite differentiability.

Our proof of Theorem \ref{MainThm}, except for a few standard reductions at the beginning, doesn't use bilinear method. Instead, as in \cite{W15}, the endgame of the proof relies on the theorem of Matom\"aki-Radziw\l\l{}-Tao. To be accurate, we use an application that Matom\"aki-Radziw\l\l{}-Tao derived from their theorem, saying that when $R$ is sufficiently large and $N$ is sufficiently large compared to $R$, for most $1\leq L\leq N$, the correlation between $\mu(n)$ and $e(\alpha n)$ on the short interval $[L, L+R]$ is small (see Proposition \ref{MRTProp}). However, unlike in \cite{W15}, the use of \cite{MRT15} in this paper is only qualitative, eliminating the need of analyticity.

The dynamical analysis in this paper is quite different from that in \cite{W15}. The main strategy in \cite{W15} was to prove that the dynamics is almost periodic at a single step length for a very long time. In contrast, to prove Theorem \ref{MainEq}, we analyse the trajectory at multiple scales and show that it has a structure that looks like a sum of independent random variables at different scales. These scales are determined by the continued fraction expansion of $\alpha$. Thanks to the mutual independence, in order to have a measurable invariant section, the total variance of these random variables have to be bounded. This allows to study the dynamics only on the first finitely many scales. After such a reduction, the dynamics resembles a linear flow when $n$ is restricted to a short interval of given length. At this stage, \cite{MRT15} can be applied.

\noindent {\bf Notations.}\begin{itemize}
 \item $\rmm_{\bT^d}=$ the Lebesgue probability measure on $\bT^d$;
 \item $e(\theta)=e^{2\pi i\theta}$;
 \item $\|\theta\|=\mathrm{dist}(\theta,\bZ)$ for $\theta\in\bR$. Remark that $\|\theta\|\ll e(\theta)-1\ll\|\theta\|$;
 \item $\Exp_F$, $\Var_F$, $\Prob_F$ respectively stand for the expectation, variance and probability of a function/event defined on a finite set $F$ with respect to the uniform probability measure on $F$.
\end{itemize}

\noindent {\bf Acknowledgments.} I am grateful to Jean Bourgain, Jon Chaika and Peter Sarnak for helpful discussions. This work was done when the author was a member at the Institute for Advanced Study, and I would like to thank the IAS for its support. The author was also supported by NSF grants DMS-1451247 and DMS-1501295.

\section{Elimination of non-resonant frequencies}

In the remainder of the paper, we will assume that $T$ is given by \eqref{DynaEq} where $h:\bT^1\to \bT^1$ is $C^\tau$ regular with $\tau\in(2,\infty)$.

Liu and Sarnak \cite{LS15} showed that if $\alpha$ is rational, then Conjecture \ref{SarnakConj} holds for $T$. So we will always assume $\alpha$ is irrational.

Since $T$ has a measurable invariant section, by a result of Furstenberg \cite{F61}*{Lemma 2.2} $h$ is homotopically trivial. Under this restriction, $h$ can be realized as a $C^\tau$ function from $\bT^1$ to $\bR$ and be written as
\begin{equation}\label{FourierEq}h(x)=\sum_{m\in\bZ}\hh(m)e(mx),\end{equation} where the convergence is uniform and the equality holds pointwise for all $x\in\bT^1$. Moreover, as $h$ is $C^\tau$, we have
\begin{equation}\label{FourierBdEq}|\hh(m)|\ll |m|^{-\tau}, \forall m\in\bZ\backslash\{0\}.\end{equation}

Take the continued fraction expansion $$\alpha=[0;a_1,a_2,\cdots]=\dfrac1{a_1+\dfrac1{a_2+\dfrac1{a_3+\cdots}}}$$ and let $\frac{p_k}{q_k}=[0;a_1,\cdots,a_{k-1}]$ be the corresponding convergents\footnote{It should be noted that our enumeration here differs from the more commonly used one in the literature, namely $\frac{p_k}{q_k}=[0;a_1,\cdots,a_k]$.}. Because $\alpha$ is irrational, the expansion is infinite. Denote $\theta_k=q_k\alpha-p_k$.

\begin{remark}\label{DioRmk}The following standard facts can be found in \cite{K97}:
\begin{enumerate}
\item $p_1=0$, $q_1=1$; $p_2=1$, $q_2=a_1$; $p_{k+1}=a_kp_k+p_{k-1}$ and  $q_{k+1}=a_kq_k+q_{k-1}$ for $k\geq 2$;
\item $p_k$ is coprime to $q_k$;
\item $\frac1{q_{k+1}+q_k}<\|q_k\alpha\|<\frac1{q_{k+1}}$;
\item If $|\alpha-\frac pq|<\frac1{2q^2}$ for $p,q\in\bZ$, then $\frac pq$ coincides with one of the $\frac {p_k}{q_k}$'s.
\item The sequence $\{\theta_k\}$ has alternating signs and $\|q_k\alpha\|=|\theta_k|$.
\end{enumerate}
\end{remark}

One consequence to Remark \ref{DioRmk}.(1) is that $q_{k+2}\geq q_{k+1}+q_k\geq 2q_k$.  In particular, $q_k$ grows exponentially: \begin{equation}\label{ExpGrowthEq}q_k\gg 2^{\frac k2}\text{ and }q_{k+j}\gg 2^{\frac j2}q_k.\end{equation}

Set \begin{equation}\label{FourierSupportEq}M=\bigcup_{k: q_{k+1}>q_k^{\frac\tau 2}}\{\pm m_kq_k: m_k=1,\cdots,a_k\}\end{equation}

What we prove next is essentially \cite{LS15}*{Lemma 4.1}, with slightly finer estimates.

\begin{lemma}\label{DioFrequencyLem} If \eqref{FourierBdEq} holds, then $$\sum_{m\notin M}\hh(m)\frac{1}{e(m\alpha)-1}e(mx)$$ converges uniformly to a continuous function $\psi(x)$. \end{lemma}

\begin{proof} If $m\notin M$, then we are in at least one of the three situations below:

(1) $m=0$. Since this involves only one frequecy, we can ignore this case in the study of convergence.

(2) For some $k$, $q_k\leq |m|<q_{k+1}$ but $q_k\nmid |m|$. Then $\|m\alpha\|\geq\frac1{2|m|}$. This is because otherwise, by Remark \ref{DioRmk}, $|m|=aq_j$ and $\|m\alpha\|=|m\alpha-ap_j|$ for some index $j\leq k$ and $a\in\bZ$. Since $q_k\nmid |m|$, $j<k$. Hence we have $$\|m\alpha\|=|a|\cdot\|q_j\alpha\|>\frac {|a|}{q_{j+1}+q_j}\geq \frac 1{2q_k}\geq \frac1{2m}$$ anyway in this case.

Therefore for any given $k$
\begin{equation}\label{DioFrequencyLemEq1}\begin{aligned}
&\sum_{\substack{q_k\leq |m|< q_{k+1}\\q_k\nmid m}}\left|\hh(m)\frac{1}{e(m\alpha)-1}e(mx)\right|\\
\ll &\sum_{\substack{q_k\leq |m|< q_{k+1}\\q_k\nmid m}} (|m|^{-\tau}\cdot |m|\cdot 1)
\ll \sum_{m=q_k}^{q_{k+1}}m^{-(\tau-1)}\\
\ll& q_k^{-(\tau-2)}-q_{k+1}^{-(\tau-2)}.
\end{aligned}\end{equation}

(3) $m=\pm m_kq_k$ where $m_k\in\{1,\cdots,a_k\}$ but $q_{k+1}\leq q_k^{\frac\tau 2}$. Since $m_k\|q_k\alpha\|\leq a_k\frac1{q_{k+1}}<\frac1{q_k}$, $\|m\alpha\|$ is given by $m_k\|q_k\alpha\|$ for $k\geq 2$. Thus, we have for all $k\geq 2$ that

\begin{equation}\label{DioFrequencyLemEq2}\begin{aligned}
&\sum_{\substack{q_k\leq |m|< q_{k+1}\\q_k| m}}\left|\hh(m)\frac{1}{e(m\alpha)-1}e(mx)\right|\\
\ll &2\sum_{m_k=1}^{a_k}\big( (m_kq_k)^{-\tau}\cdot\frac1{m_k\cdot \frac1{q_{k+1}+q_k}}\cdot 1\big)\\
\ll &\sum_{m_k=1}^{a_k}m_k^{-(\tau+1)}q_{k+1}q_k^{-\tau}
\ll\sum_{m_k=1}^{\infty}m_k^{-(\tau+1)}q_k^{-\frac\tau2}\\
\ll& q_k^{-\frac\tau2}.
\end{aligned}\end{equation}
The last inequality here is because $\tau-1>1$.

We sum both estimates \eqref{DioFrequencyLemEq1} and \eqref{DioFrequencyLemEq2} over all $k\geq 2$. Since only finitely many terms are neglected in doing this, and the estimates are independent of $x$, to prove the lemma it suffices to know that both the resulting series are convergent. This is indeed the case, respectively because $\tau-2>0$ and $\frac\tau2>1$. \end{proof}

\begin{corollary}\label{FourierReduction} Assuming \eqref{FourierBdEq}, Conjecture \ref{SarnakConj} holds for $T$ if and only it holds for the map $T_1(x,y)=(x+\alpha, y+h_1(x))$ on $\bT^2$, where $$h_1(x)=\sum_{m\in M\cup\{0\}}\hh(m)e(mx).$$\end{corollary}

In other words, in order to prove Theorem \ref{MainThm}, one may assume that $\hh$ is supported on $M\cup\{0\}$. The proof of the corollary is the same as that of \cite{W15}*{Corollary 3.3}. In fact, it suffices to notice that $T$ is continously conjugate to $T_1$ by the map $(x,y)\mapsto (x,y+\psi(x))$.

In the same spirit, one may assume that $M$ is infinite. In fact, if $M$ is finite then $$\sum_{m\neq 0}\hh(m)\frac{1}{e(m\alpha)-1}e(mx)$$ differs from the series in Lemma \ref{DioFrequencyLem} by only finitely many terms and hence also defines a continuous function $\tpsi$. And $T$ is continously conjugate to the Kronecker flow $(x,y)\mapsto (x+\alpha, y+\hh(0))$ by $(x,y)\mapsto (x,y+\tpsi(x))$. So it suffices to show \eqref{MainEq} for this linear flow, which is known by the work of Davenport \cite{D37}.

To summarize the reductions above, in order to prove Theorem \ref{MainThm}, we may assume the following:

\begin{hypothesis}\label{ReducedHypo}$T=(x+\alpha,y+h(x))$ with $\alpha$ irrational and $h\in C^\tau$ where $\tau>2$. Moreover, \begin{enumerate}                                                                                                                             \item $h$ is homotopically trivial;
\item The subset $M\subset\bZ$ defined by \eqref{FourierSupportEq} is infinite and $\hh$ is supported on $M$;
\item The map $T$ preserves a measurable section $(x,g(x))$.                                                                                                               \end{enumerate}
\end{hypothesis}

We will always work under Hypothesis \ref{ReducedHypo} hereafter.

\section{Approximation of trajectories: non-zero frequencies}

For $n\in\bZ$, write the $n$-th iterate of $T$ as \begin{equation}\label{IterateEq}T^n(x,y)=\big(x+n\alpha, y+H_n(x)\big).\end{equation}

Then for $n\geq 0$, $$H_n(x)=\sum_{l=0}^{n-1}h(x+l\alpha)=\sum_{l=0}^{n-1}
\sum_{m\in\bZ}\hh(m)e(mx)e(lm\alpha).$$
By exchanging the order of sums, we get \begin{equation}\label{HEq}\begin{aligned}H_n(x)=&\sum_{m\in\bZ}\hh(m)\big(\sum_{l=0}^{n-1}e(lm\alpha)\big)e(mx)\\
=&n\hh(0)+\sum_{m\in\bZ\backslash\{0\}}\hh(m)\frac{e(nm\alpha)-1}{e(m\alpha)-1}e(mx).\end{aligned}\end{equation}

When $n<0$, $$H_n(x)=\sum_{l=n}^{-1}-h(x+l\alpha)=\sum_{l=n}^{-1}\sum_{m\in\bZ}-\hh(m)e(mx)e(lm\alpha).$$ In this case, one can easily check that \eqref{HEq} remains true. Thus \eqref{HEq} holds for all $n\in\bZ$.

Given $x\in\bT^1$, we shall choose its unique representative from $[-\alpha,1-\alpha)$, which we denote by $x$ indifferently. To better approximate the series \eqref{HEq} we expand $n$ and $x$ in their Ostrowski numerations \cite{O22} with respect to $\alpha$, the definitions of which we now recall.

\begin{definition}\label{OstrowskiZDef} {\rm (Ostrowski numeration for integers)} Every non-negative integer $n$ can be uniquely written as a sum $\sum_{k\geq 1}n_kq_k$, where:\begin{enumerate}
\item $n_k$ is an integer, $0\leq n_k\leq a_k$. When $k=1$, $0\leq n_1\leq a_1-1$;                                                                                                                                                         \item $n_k=0$ if $n_{k+1}=a_{k+1}$;
\item $n_k=0$ for all but finitely many $k$'s.
\end{enumerate}
\end{definition}

\begin{definition}{\rm (Ostrowski numeration for real values)} Every value $x\in[-\alpha,1-\alpha)$ can be uniquely written as a convergent series $\sum_{k\geq 1}\tx_k\theta_k$, where:\begin{enumerate}
\item $\tx_k$ is an integer, $0\leq \tx_k\leq a_k$ in general, but $0\leq \tx_1\leq a_1-1$;                                                                                                                                                         \item $\tx_k=0$ if $\tx_{k+1}=a_{k+1}$;
\item $\tx_k<a_k$ infinitely often.\end{enumerate}
We denote $x_k=\tx_k\theta_k$.\end{definition}

For an introduction to Ostrowski numerations, see the survey \cite{B01}.

Given positive integers $2\leq k_-\leq k_+$, we are interested in estimating \eqref{HEq} for positive integers $n$ whose Ostrowski numerations have the form $n=\ksum n_kq_k$.

More generally, we may assume
\begin{equation}\label{nTruncateEq}n=\ksum n_kq_k,\ n_k\in\bZ,\ |n_k|\ll a_k,\end{equation}
without requiring the decomposition to be an Ostrowski numeration.

With \eqref{nTruncateEq}, we will denote partial sums by
\begin{equation}\label{PartialSumEq}\bn_k=\sum_{j=k_-}^kn_jq_j.\end{equation}

The main idea of this section is that the interaction between the component $n_kq_k$ of $n$ and the Fourier component at frequency $m_jq_j$ of $h$ really matters for the dynamics of $T$ only when $k=j$.

\begin{lemma}\label{ApproxLem1} For $k_-$, $k_+$,  $n$ as in \eqref{nTruncateEq} and any $x\in\bT^1$, $H_n(x)$ is approximated by the sum
$$H_n^{(1)}(x)=n\hh(0)+\ksum \sum_{\substack{-a_k\leq m_k\leq a_k\\m_k\neq 0}}\hh(m_kq_k)\frac{e(n_km_kq_k^2\alpha)-1}{e(m_kq_k\alpha)-1}e(m_kq_k(x+\bn_{k-1}\alpha)),$$ with an error term of order $O(2^{-\frac{k_-}4})$.
\end{lemma}
\begin{proof} By definition $H_n(x)=\ksum H_{n_kq_k}(x+\bn_{k-1}\alpha)$. Since we are assuming Hypothesis \ref{ReducedHypo}, it follows that
\begin{equation}\label{HEq2}\begin{aligned}
&H_n(x)\\
=&n\hh(0)+\ksum \sum_{j=1}^\infty\sum_{\substack{-a_j\leq m_j\leq a_j\\m_j\neq 0}}\hh(m_jq_j)\frac{e(n_kq_km_jq_j\alpha)-1}{e(m_jq_j\alpha)-1}e(m_jq_j(x+\bn_{k-1}\alpha)).
\end{aligned}\end{equation}
So $H_n(x)-H_n^{(1)}(x)$ is equal to $$\ksum \sum_{\substack{j\geq 1\\j\neq k}}\sum_{\substack{-a_j\leq m_j\leq a_j\\m_j\neq 0}}\hh(m_jq_j)\frac{e(n_kq_km_jq_j\alpha)-1}{e(m_jq_j\alpha)-1}e(m_jq_j(x+\bn_{k-1}\alpha)).$$

Therefore
\begin{equation}\label{ApproxLem1Eq1}\begin{aligned}
&|H_n(x)-H_n^{(1)}(x)|\\
\leq&\ksum \sum_{\substack{j\geq 1\\j\neq k}}\left|\sum_{\substack{-a_j\leq m_j\leq a_j\\m_j\neq 0}}\hh(m_jq_j)\frac{e(n_kq_km_jq_j\alpha)-1}{e(m_jq_j\alpha)-1}e(m_jq_j(x+\bn_{k-1}\alpha))\right|\\
\ll&\ksum \sum_{\substack{j\geq 1\\j\neq k}}\sum_{\substack{-a_j\leq m_j\leq a_j\\m_j\neq 0}}|\hh(m_jq_j)|\frac{\|n_kq_km_jq_j\alpha\|}{\|m_jq_j\alpha\|}
\end{aligned}\end{equation}

We distinguish between the cases $j<k$ and $j>k$.

When $j<k$,
\begin{equation}\label{ApproxLem1Eq2}\begin{aligned}
   &\sum_{\substack{-a_j\leq m_j\leq a_j\\m_j\neq 0}}|\hh(m_jq_j)|\frac{\|n_kq_km_jq_j\alpha\|}{\|m_jq_j\alpha\|}\\
\ll&\sum_{\substack{-a_j\leq m_j\leq a_j\\m_j\neq 0}}|m_j|^{-\tau}q_j^{-\tau}\frac{|n_k|\cdot |m_j|\cdot q_j\cdot|\theta_k|}{|m_j|\cdot|\theta_j|}\\
\ll&\sum_{\substack{-a_j\leq m_j\leq a_j\\m_j\neq 0}}|m_j|^{-\tau}q_j^{-\tau}\frac{|n_k|q_j q_{j+1}}{q_{k+1}}
\ll q_j^{-(\tau-1)}\frac{q_{j+1}}{q_k}
\end{aligned}\end{equation}
Here we used $q_{j+1}^{-1}\ll|\theta_j|\ll q_{j+1}^{-1}$, which is guaranteed by Remark \ref{DioRmk}, and that $\frac{|n_k|}{q_{k+1}}\ll\frac{a_k}{q_{k+1}}<\frac1{q_k}$.

One can further bound the estimate above using the exponential growth rate from \eqref{ExpGrowthEq}. If $1\leq j<\frac k2$, then $\eqref{ApproxLem1Eq2}\ll \frac{q_{j+1}}{q_k}\ll 2^{-\frac{k-j-1}2}$. If $\frac k2\leq j<k$, then  $\eqref{ApproxLem1Eq2}\ll q_j^{-(\tau-1)}\ll 2^{-\frac{j(\tau-1)}2}$ by \eqref{ExpGrowthEq}. It follows that
\begin{equation}\label{ApproxLem1Eq3}
\begin{aligned}
&\ksum \sum_{j=1}^{k-1}\sum_{\substack{-a_j\leq m_j\leq a_j\\m_j\neq 0}}|\hh(m_jq_j)|\frac{\|n_kq_km_jq_j\alpha\|}{\|m_jq_j\alpha\|}\\
\ll&\ksum \left(\sum_{j=1}^{\lceil \frac k2\rceil-1}2^{-\frac{k-j-1}2}+\sum_{j=\lceil \frac k2\rceil}^{k-1}2^{-\frac{j(\tau-1)}2}\right)\\
\ll&\ksum (2^{-\frac k4}+2^{-\frac{(\tau-1)k}4})\ll \ksum 2^{-\frac k4}
\ll 2^{-\frac{k_-}4}.
\end{aligned}\end{equation}

And when $j>k$,
\begin{equation}\label{ApproxLem1Eq4}\begin{aligned}
   &\sum_{\substack{-a_j\leq m_j\leq a_j\\m_j\neq 0}}|\hh(m_jq_j)|\frac{\|n_kq_km_jq_j\alpha\|}{\|m_jq_j\alpha\|}\\
\ll&\sum_{\substack{-a_j\leq m_j\leq a_j\\m_j\neq 0}}|m_j|^{-\tau}q_j^{-\tau}n_kq_k\ll q_j^{-\tau}q_{k+1}\\
\ll& q_j^{-\tau}q_j=q_j^{-(\tau-1)}.
\end{aligned}\end{equation}

Thus, once again thanks to the exponential growth of $\{q_k\}$,
\begin{equation}\label{ApproxLem1Eq5}
\begin{aligned}
&\ksum \sum_{j=k+1}^{\infty}\sum_{\substack{-a_j\leq m_j\leq a_j\\m_j\neq 0}}|\hh(m_jq_j)|\frac{\|n_kq_km_jq_j\alpha\|}{\|m_jq_j\alpha\|}\\
\ll&\ksum \sum_{j=k+1}^\infty q_j^{-(\tau-1)}
\ll\ksum  q_{k+1}^{-(\tau-1)}\\
\ll&q_{k_-+1}^{-(\tau-1)}.
\end{aligned}\end{equation}

Because $q_{k_-+1}\gg 2^{\frac {k_-}2}$ and $\tau-1>1$, \eqref{ApproxLem1Eq5} is dominated by \eqref{ApproxLem1Eq3}. By feeding both of them into \eqref{ApproxLem1Eq1}, we obtain that $|H_n(x)-H_n^{(1)}(x)|\ll 2^{-\frac{k_-}4}$.
\end{proof}

From \eqref{ApproxLem1Eq4} we deduce an approximation by truncated Fourier series, that will become useful in a later part of this paper.

\begin{lemma}\label{TruncApproxLem}Suppose $n=\sum_{k=k_-}^{k_+}n_kq_k$ where $|n_k|\ll a_k$, define a truncaction of $h(x)$ by $h^*(x)=\hh(0)+\sum_{j=1}^{k_+}\sum_{\substack{-a_k\leq m_k\leq a_k\\m_k\neq 0}}\hh(m_kq_k)e(m_kq_kx)$.
Also define $H^*_n(x)$ in the same way as $H_n(x)$ using $h^*(x)$ instead of $h(x)$, then
$$\|H^*_n(x)-H_n(x)\|\ll q_{k_++1}^{-(\tau-1)}.$$
\end{lemma}

\begin{proof} The difference between $H_n(x)$ and $H^*_n(x)$ is the sum of all terms in \eqref{HEq} involving frequencies $m_jq_j$, $j>k_+$, or
$$\ksum \sum_{j=k_++1}^\infty\sum_{\substack{-a_j\leq m_j\leq a_j\\m_j\neq 0}}\hh(m_jq_j)\frac{e(n_kq_km_jq_j\alpha)-1}{e(m_jq_j\alpha)-1}e(m_jq_j(x+\bn_{k-1}\alpha)).$$

By the second line in \eqref{ApproxLem1Eq4},
$$\begin{aligned}
\|H^*_n(x)-H_n(x)\|
\ll &\ksum\sum_{j=k_++1}^\infty q_j^{-\tau}q_{k+1}
\ll \ksum q_{k_++1}^{-\tau}q_{k+1}\\
\ll& q_{k_++1}^{-\tau}q_{k_++1}=q_{k_++1}^{-(\tau-1)}.
\end{aligned}$$
\end{proof}

The approximation $H_n^{(1)}(x)$ can be further refined by exploiting the Ostrowski numeration of $x$.

\begin{lemma}\label{ApproxLem2} For $k_-$, $k_+$, and $n$ as in \eqref{nTruncateEq}, $H_n(x)$ is approximated by the sum
$$H_n^{(2)}(x)=n\hh(0)+\ksum \sum_{\substack{-a_k\leq m_k\leq a_k\\m_k\neq 0}}\hh(m_kq_k)\frac{e((n_k+\tx_k)m_kq_k^2\alpha)-e(\tx_km_kq_k^2\alpha)}{e(m_kq_k\alpha)-1},$$ with an error term of order $O(2^{-\frac{k_-}4})$.  Here $\{\tx_k\}$ is the Ostrowski numeration of $x$.\end{lemma}

\begin{proof}The proof is in a sense symmetric to that of Lemma \ref{ApproxLem2}, by interchanging perspectives of the variables $n$ and $x$.

Write first
\begin{equation}\begin{aligned}
x+\bn_{k-1}\alpha=&\sum_{j=1}^\infty\tx_j \theta_j+\sum_{j=1}^{k-1} n_jq_j\alpha\\
\equiv& \tx_kq_k\alpha+\sum_{j=1}^{k-1}(n_j+\tx_j)q_j\alpha+\sum_{j=k+1}^\infty\tx_jq_j\alpha\  (\text{mod }1).
\end{aligned}\end{equation}
For simplicity, write $$\tx'_j=\tx_j+\mathbf 1_{j\leq k}n_j,$$ and $$z_j=\tx_kq_k\alpha+\sum_{l=1}^{\min(j,k-1)}\tx'_lq_l\alpha+\sum_{l=k+1}^j\tx'_lq_l\alpha.$$ Then $z_0=\tx_kq_k\alpha$; and $z_j\to x+\bn_{k-1}\alpha$ as $j\to\infty$ in $\bT^1$; and $z_j-z_{j-1}=\tx'_jq_j\alpha$ except when $j=k$, for which $z_k=z_{k-1}$. Also notice that $|\tx'_j|\ll 2a_j$.

Therefore we have a decomposition
\begin{equation}\label{ApproxLem2Eq1}\begin{aligned}
&e(m_kq_k(x_k+\bn_{k-1}\alpha))\\
=&e(m_kq_k\tx_kq_k\alpha)+\sum_{\substack{j\geq 1\\ j\neq k}}e(m_kq_kz_j)\big(e(m_kq_k\tx'_jq_j\alpha)-1\big).
\end{aligned}\end{equation} This series converges uniformly as with $k$ fixed, the $j$-th term is of the same order as $\|\tx_j'q_j\alpha\|\ll a_j\theta_j\ll a_jq_{j+1}^{-1}\ll q_j^{-1}$ and decays exponentially fast.

After plugging \eqref{ApproxLem2Eq1} into  the expression of $H_n^{(1)}(x)$ and comparing with $H_n^{(2)}(x)$, we see that

\begin{equation}\label{ApproxLem2Eq2}\begin{aligned}
&|H_n^{(1)}(x)-H_n^{(2)}(x)|\\
\leq &\Bigg|\ksum \sum_{\substack{-a_k\leq m_k\leq a_k\\m_k\neq 0}}\Big(\hh(m_kq_k)\frac{e(n_km_kq_k^2\alpha)-1}{e(m_kq_k\alpha)-1}\cdot\\
&\ \ \ \sum_{\substack{j\geq 1\\j\neq k}}e(m_kq_kz_j)\big(e(m_kq_k\tx'_jq_j\alpha)-1\big)\Big)\Bigg|\\
\ll&\ksum \sum_{\substack{-a_k\leq m_k\leq a_k\\m_k\neq 0}} \sum_{\substack{j\geq 1\\j\neq k}}|\hh(m_kq_k)|\cdot\Big|\frac{e(m_kq_k\tx'_jq_j\alpha)-1}{e(m_kq_k\alpha)-1}\Big|\\
\ll&\ksum \sum_{\substack{-a_k\leq m_k\leq a_k\\m_k\neq 0}} \sum_{\substack{j\geq 1\\j\neq k}}|\hh(m_kq_k)|\frac{\|m_kq_k\tx'_jq_j\alpha\|}{\|m_kq_k\alpha\|}.
\end{aligned}\end{equation}

We again distinguish the cases $j<k$ and $j>k$.

Suppose first $j<k$. Similar to \eqref{ApproxLem1Eq4} we have

\begin{equation}\label{ApproxLem2Eq3}\begin{aligned}
&\sum_{\substack{-a_k\leq m_k\leq a_k\\m_k\neq 0}}|\hh(m_kq_k)|\frac{\|m_kq_k\tx'_jq_j\alpha\|}{\|m_kq_k\alpha\|}\\
\ll&\sum_{\substack{-a_k\leq m_k\leq a_k\\m_k\neq 0}}|m_k|^{-\tau}q_k^{-\tau}\cdot|\tx'_j|q_j
\ll q_k^{-\tau}q_{j+1}.
\end{aligned}\end{equation}
Here we used that $|\tx'_j|q_j\ll a_jq_j<q_{j+1}$.

Because $q_k$ has exponential growth, it follows that

\begin{equation}\label{ApproxLem2Eq4}\begin{aligned}
&\ksum \sum_{\substack{-a_k\leq m_k\leq a_k\\m_k\neq 0}} \sum_{j=1}^{k-1}|\hh(m_kq_k)|\frac{\|m_kq_k\tx'_jq_j\alpha\|}{\|m_kq_k\alpha\|}\\
\ll&\ksum q_k^{-\tau}q_k=\ksum q_k^{-(\tau-1)}
\ll q_{k_-}^{-(\tau-1)}.
\end{aligned}\end{equation}

Assume now $j>k$. Because $|\tx_j'|\ll a_j$,  similar to \eqref{ApproxLem1Eq2} we have
\begin{equation}\label{ApproxLem2Eq5}\begin{aligned}
   &\sum_{\substack{-a_k\leq m_k\leq a_k\\m_k\neq 0}}|\hh(m_kq_k)|\frac{\|m_kq_k\tx'_jq_j\alpha\|}{\|m_kq_k\alpha\|}\\
\ll&\sum_{\substack{-a_k\leq m_k\leq a_k\\m_k\neq 0}}|m_k|^{-\tau}q_k^{-\tau}\frac{a_j|m_k|q_k\cdot|\theta_j|}{|m_k|\cdot|\theta_k|}\\
\ll&\sum_{\substack{-a_k\leq m_k\leq a_k\\m_k\neq 0}}|m_k|^{-\tau}q_k^{-\tau}\frac{a_j q_k q_{k+1}}{q_{j+1}}
\ll q_k^{-(\tau-1)}\frac{q_{k+1}}{q_j}.
\end{aligned}\end{equation}

Hence
\begin{equation}\label{ApproxLem2Eq6}\begin{aligned}
   &\ksum \sum_{\substack{-a_k\leq m_k\leq a_k\\m_k\neq 0}} \sum_{j=k+1}^\infty|\hh(m_kq_k)|\frac{\|m_kq_k\tx'_jq_j\alpha\|}{\|m_kq_k\alpha\|}\\
\ll&\ksum \sum_{j=k+1}^\infty q_k^{-(\tau-1)}q_{k+1}q_j^{-1}
\ll\ksum q_k^{-(\tau-1)}q_{k+1}q_{k+1}^{-1}\\
\ll& q_{k_-}^{-(\tau-1)}.
\end{aligned}\end{equation}

By adding \eqref{ApproxLem2Eq4} and \eqref{ApproxLem2Eq6} and comparing with \eqref{ApproxLem2Eq2}, we know that \begin{equation}\label{ApproxLem2Eq7}|H_n^{(1)}(x)-H_n^{(2)}(x)|\ll q_{k_-}^{-(\tau-1)}.\end{equation}

Because $\tau-1>1$ and $q_{k_-}\gg 2^{\frac{k_-}2}$, $q_{k_-}^{-(\tau-1)}\ll 2^{-\frac{k_-}4}$. The lemma is verified by combining \eqref{ApproxLem2Eq7} with Lemma \ref{ApproxLem1}.
\end{proof}

When $|n_k|$ is small enough, $H_n$ can be approximated directly by $n\hh(0)$.

\begin{lemma}\label{StepBdLem} For $k_-$, $k_+$, $n$ and $\{n_k\}$ as in \eqref{nTruncateEq},
$$\left|\ksum \sum_{\substack{-a_k\leq m_k\leq a_k\\m_k\neq 0}}\hh(m_kq_k)\frac{e(n_km_kq_k^2\alpha)-1}{e(m_kq_k\alpha)-1}\right|\ll \ksum |n_k|q_k^{-(\tau-1)}.$$
and for all $x\in\bT^1$,
 $$|H_n(x)-n\hh(0)|\ll \ksum |n_k|q_k^{-(\tau-1)}+2^{-\frac{k_-}4}.$$
\end{lemma}

\begin{proof} The first quantity is bounded by
$$\ksum \sum_{\substack{-a_k\leq m_k\leq a_k\\m_k\neq 0}}|m_k|^{-\tau}q_k^{-\tau}\cdot |n_k|q_k
\ll  \ksum |n_k|q_k^{-(\tau-1)}.$$ For the second inquality, it suffices to note $$|H_n^{(1)}-n\hh(0)|\leq \ksum \sum_{\substack{-a_k\leq m_k\leq a_k\\m_k\neq 0}}\Big|\hh(m_kq_k)\frac{e(n_km_kq_k^2\alpha)-1}{e(m_kq_k\alpha)-1}\Big|$$ and combine with Lemma \ref{ApproxLem1}.
\end{proof}

\section{Diophantine relation between rotation numbers}

The map $T$ has horizontal rotation number $x$ and vertical rotation number $\hh(0)$. The main result of this part, Proposition \ref{LuzinProp}, shows that these two numbers have similar Diophantine profile under our assumption that $T$ has a measurable invariant section.

Assume for a measurable function $g:\bT^1\mapsto\bT^1$,  \begin{equation}\label{InvSectEq}g(x+\alpha)-g(x)=h(x), \text{ a.e. }x\in\bT^1.\end{equation}

\begin{lemma}\label{InvSectLem} There is $\xi_1\in\bZ$, such that for the function $e_{(\xi_1,1)}(x,y)=e(\xi_1x+y)$, for $\rmm_{\bT^1}$-a.e. $x$ and all $y$,
$$\Exp_{n=1}^N(e_{(\xi_1,1)}\circ T^n)(x,y)\not\to 0,\text{ as }N\to\infty.$$\end{lemma}
\begin{proof} Notice that measurable function $f(x,y)=e(-g(x)+y)$ satisfies $|f(x,y)|=1$ and $f\circ T=f$ for almost every $x$ and all $y$. Decompose $e(g(x))$ as $\sum_{\xi_1\in\bZ}c_{\xi_1}e(\xi_1x)$ in $L^2$. Then $f=\sum_{\xi_1\in\bZ}c_{\xi_1}e_{(\xi_1,1)}$ in $L^2(\rmm_{\bT^2})$.

It follows that for at least one $\xi_1$, $\Exp_{n=1}^Ne_{(\xi_1,1)}\circ T^n$ does not converge to $0$ in $L^2$. Otherwise, $\Exp_{n=1}^Nf\circ T^n$ would converge to $0$ in $L^2$ as well, contradicting the fact that it is always equal to $f$.

Remark that for the same $x$, and different $y$, $y'$, $$\Exp_{n=1}^N(e_{(\xi_1,1)}\circ T^n)(x,y)=e(y-y')\cdot \Exp_{n=1}^N(e_{(\xi_1,1)}\circ T^n)(x,y').$$ In other words, given $x$, whether the ergodic average converges to $0$ is independent of $y$. Denote by $\Omega$ the set of $x$ for which the averages converge to $0$. Then $\Omega$ is invariant under $x\mapsto x+\alpha$, and have Lebesgue measure $0$ or $1$ by ergodicity. Assume $\rmm_{\bT^1}(\Omega)=1$, then $\Exp_{n=1}^Ne_{(\xi_1,1)}\circ T^n$ converges to $0$ pointwisely, and hence also in $L^2$ by dominated convergence theorem. This contradicts the choice of $x$. Hence $\rmm_{\bT^1}(\Omega)=0$. The lemma is proved.
\end{proof}

\begin{proposition}\label{LuzinProp} For all $\delta>0$, there exists $k_0=k_0(\delta)\in\bN$ with the following property:

If $k_0\leq k_-\leq k_+$, $n$ is as in \eqref{nTruncateEq}, and $|n_k|\leq \min(q_k^{\frac{\tau-1}2},a_k)$ for all $k_-\leq k\leq k_+$, then  $\|n\hh(0)\|<\delta$.\end{proposition}

The proof follows a similar approach as \cite{W15}*{Lemma 4.1}.

\begin{proof} By Luzin's theorem, there is a compact subset $\Omega\subset\bT^1$ with $\rmm_{\bT^1}(\Omega)>\frac12$, on which $g$ is continuous. There is an $\epsilon>0$, such that if $x,x'\in\Omega$ satisfies $\|x-x'\|<\epsilon$, then $\|g(x)-g(x')\|<\frac\delta{2}$.

We claim that when $k_-$ is sufficiently large, $\|n\alpha\|<\epsilon$. This is because
$$\begin{aligned}\|n\alpha\|=&\left\|\sum_{k=k-}^{k_+}n_kq_k\alpha\right\|\leq\ksum |n_k|\cdot|\theta_k|<\ksum a_kq_{k+1}^{-1}<\ksum q_k^{-1}\\
\ll& q_k^{-1}. \end{aligned}$$

Now assume $k_-$ is large enough and thus $\|n\alpha\|<\epsilon$. The set $\{x\in\bT^1: x,x+n\alpha\in\Omega\}$ at least has Lebesgue measure $1-2\rmm_{\bT^1}(\Omega^c)>0$, and is thus non-empty. Fix a point $x$ from this set, then $\|g(x+n\alpha)-g(x)\|<\frac\delta{2}$ by the choice of $\epsilon$.

Since $g(x+n\alpha)=g(x)+H_n(x)$ is true almost everywhere, we can assume it is satisfied by the chosen $x$. Therefore, $\|H_n(x)\|<\frac\delta{2}$.

On the other hand, by Lemma \ref{ApproxLem2} and Lemma \ref{StepBdLem},
$$\begin{aligned}&|H_n(x)-n\hh(0)|\\
\ll&\left(\ksum q_k^{\frac{\tau-1}2}q_k^{-(\tau-1)}+2^{-\frac{k_-}4}
\right)=\left(\ksum q_k^{-\frac{\tau-1}2}+2^{-\frac{k_-}4}\right)\\
\ll &(q_{k_-}^{-\frac{\tau-1}2}+2^{-\frac{k_-}4})
\end{aligned}
$$ is less than $\frac\delta{2}$ when $k_-$ is sufficiently large.

By adding the two estimates above, we see that, given $\delta$, $$\|n\hh(x)\|\leq \|H_n(x)\|+|H_n(x)-n\hh(0)|<\frac\delta{2}+\frac\delta{2}= \delta$$ when $k_-$ is sufficiently large. This completes the proof.
\end{proof}

\section{Trajectories as sums of random variables}

Lemma \ref{ApproxLem2} provides an approximation $H_n^{(2)}$ to $H_n(x)$. Next, we examine $H_n^{(2)}$ more carefully and show that it is approximately the sum of a sequence of independent random variables:

\begin{equation}\label{ApproxLem2RewriteEq}H_n^{(2)}(x)=\ksum \big(\tphi_k(n_k+\tx_k)-\tphi_k(\tx_k)\big),\end{equation} where
\begin{equation}\tphi_k(l)=lq_k\hh(0)+\sum_{\substack{-a_k\leq m_k\leq a_k\\m_k\neq 0}}\hh(m_kq_k)\frac{e(lm_kq_k^2\alpha)-1}{e(m_kq_k\alpha)-1}.\end{equation}

\begin{lemma}\label{IncrePerLem}For all $\delta>0$, there exists  $k_1\in\bN$ with the following property:

Suppose $k_1\leq k_-\leq k_+$; $n=\ksum {n_k}q_k$, $n'=\ksum {n'_k}q_k$ where $n_k$, $n'_k$ are integers of absolute value bounded by $4a_k$, and $n_k\equiv n'_k\text{\rm (mod }a_k)$ for all $k_-\leq k\leq k$, then $$\left\|\ksum \big(\tphi_k(n_k)-\tphi_k(n'_k)\big)\right\|\leq\frac{3\delta}{8} . $$\end{lemma}

\begin{proof}Let $n_k-n'_k=b_ka_k$ with $|b_k|\leq 8$. Then \begin{equation}\label{IncrePerLemEq1}(n_k-n'_k)=b_ka_kq_k=b_kq_{k+1}-b_kq_{k-1}.\end{equation}
Decompose
\begin{equation}\label{IncrePerLemEq2}\begin{aligned}
\ksum &\big(\tphi_k(n)-\tphi_k(n')
= \ksum (n_k-n'_k)q_k\hh(0)+\\
&\ksum \sum_{\substack{-a_k\leq m_k\leq a_k\\m_k\neq 0}}\hh(m_kq_k)\frac{e((n_k-n'_k)m_kq_k^2\alpha)-1}{e(m_kq_k\alpha)-1}e(n'_km_kq_k^2\alpha).
\end{aligned}\end{equation}

The absolute value of the second term in \eqref{IncrePerLemEq2} is bounded by \begin{equation}\label{IncrePerLemEq3}\begin{aligned}
&\ksum \sum_{\substack{-a_k\leq m_k\leq a_k\\m_k\neq 0}}\left|\hh(m_kq_k)\frac{e((n_k-n'_k)m_kq_k^2\alpha)-1}
{e(m_kq_k\alpha)-1}\right|\\
=&\ksum \sum_{\substack{-a_k\leq m_k\leq a_k\\m_k\neq 0}}\left|\hh(m_kq_k)\frac{e((b_kq_{k+1}-b_kq_{k-1})m_kq_k\alpha)-1}
{e(m_kq_k\alpha)-1}\right|\\
\leq & \sum_{k=k_-+1}^{k_++1}\sum_{\substack{-a_{k-1}\leq m_{k-1}\leq a_{k-1}\\m_{k-1}\neq 0}}\left|\hh(m_{k-1}q_{k-1})\frac{e(b_{k-1}q_km_{k-1}q_{k-1}\alpha)-1}
{e(m_{k-1}q_{k-1}\alpha)-1}\right|\\
&+\sum_{k=k_--1}^{k_+-1}\sum_{\substack{-a_{k+1}\leq m_{k+1}\leq a_{k+1}\\m_{k+1}\neq 0}}\left|\hh(m_{k+1}q_{k+1})\frac{e(b_{k+1}q_km_{k+1}q_{k+1}\alpha)-1}
{e(m_{k+1}q_{k+1}\alpha)-1}\right|.\\
\end{aligned}\end{equation}

Here in the last inequality we used the basic fact that $|e(s_1+s_2)-1|\leq |e(s_1)-1|+|e(s_2)-1|.$

Because $|b_{k-1}|,|b_{k+1}|\leq 8\ll a_k$, both terms in \eqref{IncrePerLemEq3} are bounded by the estimate \eqref{ApproxLem1Eq1}, respectively after replacing the interval of indices $[k_-,k_+]$ with $[k_-+1,k_++1]$ and with $[k_--1,k_--1]$. From the bound on \eqref{ApproxLem1Eq1} in the proof of Lemma \ref{ApproxLem1},  we get

\begin{equation}\label{IncrePerLemEq4} \eqref{IncrePerLemEq3}\ll 2^{-\frac{k_--1}4}+2^{-\frac{k_+-1}4}\ll 2^{-\frac{k_-}4}<\frac{\delta}{8},\end{equation} when $k_1$ is sufficiently large.

We now focus on the first term from the right hand side in \eqref{IncrePerLemEq2}, which is equal to
\begin{equation}\label{IncrePerLemEq5}\begin{aligned}&\ksum (b_kq_{k+1}-b_kq_{k-1})\hh(0)\\
=&\Big(\sum_{k=k_--1}^{k_+-1}b_{k-1}q_k\Big)\hh(0)
+\Big(\sum_{k=k_--1}^{k_+-1}b_{k-1}q_k\Big)\hh(0).
\end{aligned}\end{equation}

If $k_1$ is large enough, then Proposition \ref{LuzinProp} applies to both terms since $k_--1\geq k_1$. In consequence,
\begin{equation}\label{IncrePerLemEq6}\|\eqref{IncrePerLemEq5}\|\leq\frac{\delta}{8}+\frac{\delta}{8}=\frac{\delta}{4}.\end{equation}

The lemma follows by adding \eqref{IncrePerLemEq4} to \eqref{IncrePerLemEq6}.
\end{proof}

Lemma \ref{IncrePerLem} says the projection of the function $\tphi_k$ to $\bR\mapsto \bZ$ is almost periodic with period $a_k$.

\begin{definition}\label{PerIncDef}Define a function $\phi_k:\bZ/a_k\bZ\mapsto\bT^1$ by
$$\phi_k(l\text{\rm\ mod }a_k)=\tphi_k(l)\text{\rm\ mod }\bZ, \text{ for }l=0,\cdots,a_k-1.$$\end{definition}
We shall identify $\phi_k$ with an $a_k$-periodic function on $\bZ$ without further notice. No confusion should arise from doing so.

Then Lemma \ref{IncrePerLem} actually asserts that

\begin{corollary}\label{IncrePerCor} In the settings of Lemma \ref{IncrePerLem}, $\Big\|\ksum \big(\tphi_k(n)-\phi_k(n)\big)\Big\|\leq\frac{3\delta}{8}$.\end{corollary}

We will view the digits $n_k$ from the Ostrowski numeration of an integer as random variables, and $\phi_k(n_k)$ another sequence of random variables decided by $n_k$. It will be shown in the next section that these variables are approximately independent of each other when $k$ varies.

\begin{proposition}\label{RandomSumProp}Assume Hypothesis \ref{ReducedHypo} holds. For all $\delta>0$, there exists $k_1=k_1(\delta)\geq 10$, such that:

If $k_1\leq k_-<k_+$, then for:
\begin{itemize}
 \item all integer sequences $\{\tx_k\}$ such that $|\tx_k|\leq 2a_k$ and $\sum_{k=1}^\infty\tx_k\theta_k$ converges to a real value $x$;
 \item all integers $n=\ksum n_kq_k$ where $|n_k|\leq a_k$,
\end{itemize}
we have
$$\left\|H_n(x)-\ksum \big(\phi_k(n_k+\tx_k)-\phi_k(\tx_k)\big)\right\|<\delta.$$
\end{proposition}

\begin{proof}By \eqref{ApproxLem2RewriteEq}, the left hand sided is at most
\begin{equation}\begin{aligned}
&\left\|H_n(x)-H_n^{(2)}(x)\right\|
+\left\|\ksum \big(\tphi_k(n_k+\tx_k)-\phi_k(n_k+\tx_k)\big)\right\|\\
&\hskip1cm +\left\|\ksum \big(\tphi_k(\tx_k)-\phi_k(\tx_k)\big)\right\|
\end{aligned}\end{equation}
The first term is bounded by $O(2^{-\frac{k_-}4})$ by Lemma \ref{ApproxLem2}. Because $|\tx_k|\leq 2a_k$ and $|n_k+\tx_k|\leq 3a_k$, the second term and third terms are bounded by $\frac{3\delta}{8}$ by Corollary \ref{IncrePerCor}. The proposition follows, as for sufficiently large $k_1$, one can make the $O(2^{-\frac{k_-}4})$ less than $\frac{\delta}{4}$.
\end{proof}

\section{Independence between digits in Ostrowski numerations}

In this section, we take a brief detour to show that the digits $n_k$, for indices $k$ at which $a_k$ is large, are distributed in an almost independent way for randomly chosen $n$.

\begin{definition}\label{OstrowIntervalDef} Let $I_{k_-}^{k_+}$ denote the set of integers $n\geq 0$ whose Ostrowski numeration $\{n_k\}$ is supported on $k_-\leq k\leq k_+$.
\end{definition}

Remark that
\begin{equation}I_1^k=\{0,\cdots,q_k-1\}.\end{equation}

\begin{lemma}\label{IndepLem1} Without loss of generality we can assume $k_-\geq 2$. For all pairs $k_-\leq k_+$ and for all $\psi:\bN\mapsto\bC$ with absolute value bounded by $1$,
$$\left|\Exp_{n\in I_{k_-}^{k_+}}\psi(n_{k_-})-\Exp_{l=0}^{a_{k_-}-1}\psi(l)\right|\ll\frac1{a_{k_-}}.$$
The implied constant is absolute.\end{lemma}

\begin{proof}The lemma is obvious when $k_-=k_+$, in which case $n_{k_-}$ is uniformly distributed on $\{0,\cdots,a_{k_-}\}$. We shall assume $k_-<k_+$ below.

Write $J_l=\{n\in I_{k_-}^{k_+}:n_{k_-}=l\}$ for $l=0,\cdots,a_{k_-}$. And further decompose $J_0$ as the union of $J'_0=\{n\in J_0: n_{k_-+1}\neq a_{k_-+1}\}$ and $J''_0=\{n\in J_0: n_{k_-+1}=a_{k_-+1}\}$.  Recall that if a finitely supported sequence $\{n_k\}$ is known to satisfy $0\leq n_k\leq a_k$ and $1\leq n_1\leq a_1-1$, the only restriction for it to be an Ostrowski numeration is that $n_k=0$ if $n_{k+1}=a_{k+1}$. This fact yields a few observations:
\begin{enumerate}
 \item If $1\leq l,l'\leq a_k$, then the operation of replacing $n_{k_-}$ with $l'$ in the Ostrowski numeration is a bijection from $J_l$ to $J_{l'}$;
 \item If $1\leq l\leq a_k$, then the operation of replacing $n_{k_-}$ with $l$ in the Ostrowski numeration is an injection from $J'_0$ to $J_l$.
 \item If $1\leq l\leq a_k$,  then the operation of replacing $n_{k_-}$ with $l$ and $n_{k_-+1}$ with $0$ in the Ostrowski numeration is an injection from $J''_0$ to $J_l$.
\end{enumerate}

So $|J'_0|,|J''_0|\leq |J_1|=|J_2|=\cdots=|J_{a_{k_-}}|$. Since $I_{k_-}^{k_+}$ is the disjoint union of $J'_0, J''_0, J_1, \cdots,J_{a_{k_-}}$, the lemma follows.
\end{proof}

For $\bfr=(r_0,\cdots, r_T)\in \prod_{t=0}^TI_{k_t}^{k_{t+1}-1}$, denote the concatenation of the Ostrowski numerations of the $r_t$'s by
$$\tilde\bfr=\Big((r_0)_{k_0},\cdots, (r_0)_ {k_2-1}, (r_1)_{k_2},\cdots,(r_1)_{k_3-1},\cdots, (r_T)_{k_T}, \cdots,(r_T)_{k_{T+1}-1}\Big).$$

\begin{lemma}\label{IndepLem2} Suppose $k_-=k_0<k_2<\cdots<k_{T+1}=k_++1$. If $\psi:\bN^{k_++1-k_-}\mapsto\bC$ is of absolute value bounded by $1$, then
$$\left|\Exp_{\bfr\in\prod_{t=0}^TI_{k_t}^{k_{t+1}-1}}\psi(\tilde\bfr)
-\Exp_{n\in I_{k_-}^{k_+}}\psi(n_{k_-},\cdots, n_{k_+})\right|\ll\sum_{t=1}^T\frac1{a_{k_t}}.$$
The implied constant is absolute.
\end{lemma}

\begin{proof}Define a map $S:\prod_{t=0}^TI_{k_t}^{k_{t+1}-1}\mapsto\bN$ and by $S(r_0,\cdots,r_T)=\sum_{t=0}^Tr_t$. Also define $P=(P_0,\cdots,P_T):\bN\mapsto \prod_{t=0}^TI_{k_t}^{k_{t+1}-1}$ by $P_t(n)=\sum_{k=k_t}^{k_{t+1}-1}n_kq_k$. Then the restriction of $S\circ P$ to $I_{k_-}^{k_+}$ is the identity map. Therefore, to prove the proposition, it suffices to show that the complement \begin{equation}\label{IndepLem2Eq1}\Big(\prod_{t=0}^TI_{k_t}^{k_{t+1}-1}\Big)\backslash P(I_{k_-}^{k_+})\end{equation} at most occupies an  $O(\sum_{t=1}^T\frac1{a_{k_t}})$-portion in $\prod_{t=0}^TI_{k_t}^{k_{t+1}-1}$.

Remark that $(r_0,\cdots,r_T)\in \prod_{t=0}^TI_{k_t}^{k_{t+1}-1}$ is in the set \eqref{IndepLem2Eq1} if and only if the concatenation $\tilde\bfr$ fails to be an Ostrowski numeration. This can happen only if for some $1\leq t\leq T$, $(r_{t-1})_{k_t-1}\neq 0$ but $(r_t)_{k_t}=a_{k_t}$. Hence
\begin{equation}\label{IndepLem2Eq2}\eqref{IndepLem2Eq1}\subset\bigcup_{t=1}^T\Big(J_{k_t}^{k_{t+1}-1}\times \prod_{\substack{0\leq t'\leq T\\t'\neq T}}I_{k_{t'}}^{k_{t'+1}-1}\Big),\end{equation}
where
$$J_{k_t}^{k_{t+1}-1}=\left\{r_t\in I_{k_t}^{k_{t+1}-1}: (r_t)_{k_t}\neq a_{k_t}\right\}.$$

Hence to establish the claim, one only needs to show for all $1\leq t\leq T$ that \begin{equation}
\frac{\Big|I_{k_t}^{k_{t+1}-1}\backslash J_{k_t}^{k_{t+1}-1}\Big|}{\Big|I_{k_t}^{k_{t+1}-1}\Big|}\ll\frac1{a_{k_t}}.\end{equation} This is guaranteed by the previous lemma.
\end{proof}

\begin{proposition}\label{IndepProp} For all sequences $k_-\leq k_2<\cdots <k_T\leq k_+$, and all functions $\psi:\bN^T\mapsto\bC$ with $|\psi|\leq 1$,
$$\left|\Exp_{\substack{0\leq l_t\leq a_{k_t}-1\\t=1,\cdots,T}}\psi(l_1,\cdots,l_T)-\Exp_{n\in I_{k_-}^{k_+}}\psi(n_{k_2},\cdots,n_{k_T})\right|\ll\sum_{t=1}^T\frac1{a_{k_t}},$$ where $\{n_k\}$ is the Ostrowski numeration of $n$. The implied constant is absolute.
\end{proposition}

\begin{proof} By Lemma \ref{IndepLem2}, \begin{equation}\label{IndepPropEq1}\left|\Exp_{\substack{r_t\in I_{k_t}^{k_{t+1}-1}\\t=1,\cdots,T}}\psi\big((r_1)_{k_2},\cdots,(r_T)_{k_T})-\Exp_{n\in I_{k_-}^{k_+}}\psi(n_{k_2},\cdots,n_{k_T})\right|\ll\sum_{t=1}^T\frac1{a_{k_t}}.\end{equation}

For each given $s$, by Lemma \ref{IndepLem1},
\begin{equation}\label{IndepPropEq2}
\begin{aligned}
&\Bigg|\Exp_{\substack{r_t\in I_{k_t}^{k_{t+1}-1}\\t=0,\cdots,s}}\Exp_{\substack{0\leq l_t\leq a_{k_t}-1\\t=s+1,\cdots,T}}\psi\big((r_1)_{k_2},\cdots,(r_{s-1})_{k_{s-1}},(r_s)_{k_s},l_{s+1},\cdots, l_T)\\
&\hskip1cm -\Exp_{\substack{r_t\in I_{k_t}^{k_{t+1}-1}\\t=0,\cdots,s-1}}\Exp_{\substack{0\leq l_t\leq a_{k_t}-1\\t=s,\cdots,T}}\psi\big((r_1)_{k_2},\cdots,(r_{s-1})_{k_{s-1}},l_s,l_{s+1},\cdots, l_T)\Bigg|\\
\ll&\frac 1{a_{k_s}}.\end{aligned}\end{equation}

The proposition is proved by summing \eqref{IndepPropEq2} over $s=1,\cdots, T$ and adding \eqref{IndepPropEq1}.
\end{proof}

\section{Invariant section and bounded variance}

Recall that Hypothesis \ref{ReducedHypo} is assumed throughout, and $T$ has a measurable invariant section. We will see that in consequence, the $L^2$-variances of the variables $\phi_k:\bZ/a_k\bZ\mapsto\bT^1$ are summable over $k$.

\begin{definition}Denote by $ B_{k_-}^{k_+}$  the set of finitely supported integer sequences $\bfn=\{n_k\}$ such that $0\leq n_k\leq a_k-1$ if $k_-\leq k\leq k_+$ and $q_{k+1}>q_k^{\frac\tau2}$, and $n_k=0$ otherwise.\end{definition}

The symbol $[\bfn]$ will denote $n=\sum_kn_kq_k$ for a finitely supported sequence $\bfn=\{n_k\}$.

\begin{lemma}\label{ProdDecayLem} For $k_-$ sufficiently large, the infinite product $$\prod_{\substack{k\geq k_-\\ q_{k+1}>q_k^{\frac\tau 2}}}\left|\Exp_{l=0}^{a_k-1} e\big(\phi_k(l)\big)\right|$$ converges to a non-zero value. \end{lemma}
\begin{proof} Fix an $x\in[-\alpha,1-\alpha)$ and $\xi_1$ for which Lemma \ref{InvSectLem} holds. In particular, the conclusion holds for the point $(x,0)$, i.e. the ergodic average $\Exp_{n=1}^N(e_{(\xi_1,1)}\circ T^n)(x,0)$ does not converge to $0$.
To be precise, for some $c>0$, there is a subsequence $N_i\to\infty$, such that
\begin{equation}\label{ErgAvgEq}\left|\Exp_{n=1}^{N_i}(e_{(\xi_1,1)}\circ T^n)(x,0)\right|>c.\end{equation}

For a choice of $k'_0=k'_0(c)\in\bN$ that will be specified later and $k'_0\leq k_-\leq k_+$, choose $N=N_i$ from the subsequence above, such that

\begin{equation}\label{LargeNEq}\frac{q_{k_++1}}N<\frac{c}{8}.\end{equation}

Use ensembles of the form $\{L+[\bfn]:\ \bfn\in B_{k_-}^{k_+}\}$, $L=1,\cdots, N$ to cover $\{1,\cdots,N\}$. Here elements in the ensemble are counted with multiplicity.

Each $1\leq l\leq N$ are covered exactly $| B_{k_-}^{k_+}|$ times, except for those numbers of distance less than or equal to $\max_{\bfn\in B_{k_-}^{k_+}}[\bfn]$ from either $1$ or $N$. The number of such exceptions is bounded by
$$2\sum_{k=1}^{k_+}(a_k-1)q_k=2\sum_{k=1}^{k_+}(q_{k+1}-q_{k-1}-q_k)\leq 2q_{k_++1}.$$

Therefore $\Exp_{n=1}^N\big(e_{(\xi_1,1)}\circ T^n)(x,0)$ is approximated by the decomposed average $\Exp_{L=1}^N\Exp_{\bfn\in  B_{k_-}^{k_+}}\big(e_{(\xi_1,1)}\circ T^{L+[\bfn]}\big)(x,0)$ up to an error term that is at most $\frac{2q_{k_++1}}N<\frac c{4}$. Hence by \eqref{ErgAvgEq},
$$\begin{aligned}
\Exp_{L=1}^N\left|\Exp_{\bfn\in  B_{k_-}^{k_+}}\Big(e_{(\xi_1,1)}\circ T^{L+[\bfn]}\Big)(x,0)\right|
\geq\left|\Exp_{n=1}^{N_i}(e_{(\xi_1,1)}\circ T^n)(x,0)\right|-\frac c{4}
>\frac{3c}{4}.\end{aligned}$$

In particular, there is some $L\in\{1,\cdots,N\}$ for which \begin{equation}\label{DecomposedAvgEq}\left|\Exp_{\bfn\in  B_{k_-}^{k_+}}\Big(e_{(\xi_1,1)}\circ T^{L+[\bfn]}\Big)(x,0)\right|>\frac{3c}{4}.\end{equation}

Write $T^L(x,0)=(x+L\alpha,y)$. For $\bfn\in B_{k_-}^{k_+}$,
$$\begin{aligned}(e_{(\xi_1,1)}\circ T^{L+[\bfn]})(x,0)=&e_{(\xi_1,1)}(x+(L+[\bfn])\alpha,y+H_{[\bfn]}(x+L\alpha))\\
=&e(\xi_1(x+L\alpha)+y)e(\xi_1[\bfn]\alpha+H_{[\bfn]}(x+L\alpha)).\end{aligned}$$ So \eqref{DecomposedAvgEq} is equivalent to
\begin{equation}\label{BlockAvgEq}\left|\Exp_{\bfn\in  B_{k_-}^{k_+}}e(\xi_1[\bfn]\alpha+H_{[\bfn]}(x+L\alpha))\right|>\frac{3c}{4}.\end{equation}

We will choose $k'_0$ to be greater than $k_0(\frac{c}{16\pi})$ where $k_0$ is as in Proposition \ref{RandomSumProp}. Then with $\{\tx'_k\}$ denoting the Ostrowski numeration of $x+L\alpha$,
\begin{equation}\label{BlockAvgEq1}\left\|H_{[\bfn]}(x+L\alpha)-\ksum\big(\phi_k(n_k+\tx'_k)-\phi_k(\tx'_k)\big)\right\|<\frac{c}{16\pi}.\end{equation}

In addition, for $\xi_1$ is fixed,
$$\|\xi_1[\bfn]\alpha\|=\left\|\xi_1\ksum n_kq_k\alpha\right\|\leq\xi_1\ksum a_k|\theta_k|
<\xi_1\ksum q_k^{-1}\ll q_{k_-}^{-1}.$$ So for sufficiently large $k'_0$, \begin{equation}\label{BlockAvgEq2}\|\xi_1[\bfn]\alpha\|<\frac{c}{16\pi}.\end{equation}

Using that $s\mapsto e(s)$ is $2\pi$-Lipschitz, we deduce from \eqref{BlockAvgEq}, \eqref{BlockAvgEq1}, \eqref{BlockAvgEq2} that
$$\left|\Exp_{\bfn\in  B_{k_-}^{k_+}}e\Big(\ksum\big(\phi_k(n_k+\tx'_k)-\phi_k(\tx'_k)\big)\Big)\right|>\frac{3c}{4}-2\pi(\frac{c}{16\pi}+\frac{c}{16\pi})=\frac{c}{2}.$$

By construction of $ B_{k-}^{k_+}$, the left hand side is equal to
$$\begin{aligned}&\kqprod\left|\Exp_{l=0}^{a_k-1} e\Big(\big(\phi_k(l+\tx'_k)-\phi_k(\tx'_k)\big)\Big)\right|\\
=&\kqprod\left|\Exp_{l=0}^{a_k-1} e\big(\phi_k(l+\tx'_k)\big)\right|\\
=&\kqprod\left|\Exp_{l=0}^{a_k-1} e\big(\phi_k(l)\big)\right|,\end{aligned}$$
where the last step is because $\phi_k$ is defined on $\bZ/a_k\bZ$. Hence \begin{equation}\kqprod\left|\Exp_{l=0}^{a_k-1} e\big(\phi_k(l)\big)\right|>\frac c2.\end{equation} This is true for all $k'_0\leq k_-\leq k_+$.

Given that $\frac c2>0$ and $\left|\Exp_{l=0}^{a_k-1} e\big(\phi_k(l)\big)\right|\leq 1$, the lemma follows. \end{proof}

\begin{corollary}\label{EnsembleVarCor}For all $\delta>0$, there exists $k'_1=k'_1(\delta)\in\bN$ such that for all $k_+\geq k_-\geq k'_1$, there exists a constant $z$ with $|z|\leq 1$, such that
$$\Prob_{\bfn\in  B_{k_-}^{k_+}}\left(\left|e\Big(\ksum \phi_k(n_k)\Big)-z\right|\geq\delta\right)<\delta.$$\end{corollary}
\begin{proof} It follows directly from Lemma \ref{ProdDecayLem} that, when $k_-$ is sufficiently large, for all $k_+\geq k_-$, $$\Exp_{\bfn\in B_{k_-}^{k_+}}\left|e\Big(\ksum \phi_k(n_k)\Big)\right|>\sqrt{1-\delta^3}.$$

Recall the variance for a random variable $X$ is
$$\Var(X)=\bE|X|^2-|\bE X|^2=\bE|X-\bE X|^2.$$
Since the random variable $X=e\Big(\ksum \phi_k(n_k)\Big)$ with respect to the uniform probability on $ B_{k_-}^{k_+}$ always has absolute value $1$, its variance is less than $\delta^3$. Thus
\begin{equation}\label{SmallVarPropEq2}\Prob(|X-\bE X|\geq \delta)\leq\frac{\bE|X-\bE X|^2}{\delta^2}<\delta.\end{equation} This proves the Proposition with $z=\bE X$, which satisfies $|\bE X|\leq 1$.
\end{proof}

The next step is to replace the ensemble $B_{k_-}^{k_+}$ with the natural subset $I_{k-}^{k_+}$ of $\bN$ without heavily distorting the probability distribution.

\begin{proposition}\label{IntervalVarProp}For all $\delta>0$, there exists $k_2=k_2(\delta)\in\bN$ such that for all $k_3\leq k_-\leq k_+$, there exists a constant $z$ with $|z|\leq 1$, such that
$$\Prob_{n\in I_{k_-}^{k_+}}\Big(\Big|e\Big(\kqsum \phi_k(n_k)\Big)-z\Big|\geq \delta\Big)<\delta.$$
\end{proposition}
\begin{proof} By Corollary \ref{EnsembleVarCor} and Proposition \ref{IndepProp},  for sufficiently large $k_2$ and all $k_+\geq k_-\geq k_2$, there exists $z$ with $|z|\leq 1$ such that
\begin{equation}\label{IntervalVarPropEq1}\Prob_{n\in I_{k_-}^{k_+}}\Big(\Big|e\Big(\kqsum \phi_k(n_k)\Big)-z\Big|\geq \frac\delta4\Big)<\frac\delta4+O(\kqsum\frac1{a_k}).\end{equation}

We can choose the lower bound $k_2$ such that the term $O(\kqsum\frac1{a_k})$ is less than $\frac\delta4$ for all  $k_+\geq k_-\geq k_2$. This is because for the indices $k$ involved,
\begin{equation}\label{akGrowthEq}a_k=\frac{q_{k+1}-q_{k-1}}{q_k}\geq\frac{q_{k+1}}{q_k}-1\gg q_k^{\frac\tau2-1},\end{equation} which implies $a_k$ has exponential growth with respect to $k$. Therefore \eqref{IntervalVarPropEq1} is bounded by $\frac\delta2$.

When $q_{k+1}\leq q_k^{\frac\tau2}$, by Hypothesis \ref{ReducedHypo} and the definition of $\phi_k$, $\phi_k(n_k)=\tphi_k(n_k)=n_kq_k\hh(0)$, thus $$\sum_{\substack{k_-\leq k\leq k_+\\q_{k+1}\leq q_k^{\frac\tau2}}}\phi_k(n_k)=\Big(\sum_{\substack{k_-\leq k\leq k_+\\q_{k+1}\leq q_k^{\frac\tau2}}}n_kq_k\Big)\hh(0).$$

Note $0\leq n_k\leq a_k\leq q_{k+1}q_k^{-1}\leq q_k^{\frac\tau2-1}$ if $q_{k+1}\leq q_k^{\frac\tau2}$. So by Proposition \ref{LuzinProp}, if $k_-$ is sufficiently large, then for all $n\in I_{k_-}^{k_+}$,

$$\Big|\sum_{\substack{k_-\leq k\leq k_+\\q_{k+1}\leq q_k^{\frac\tau2}}}\phi_k(n_k)\Big|<\frac\delta{4\pi},$$ which yields that
\begin{equation}\Big|e\Big(\kqsum \phi_k(n_k)\Big)-e\Big(\ksum \phi_k(n_k)\Big)\Big|<2\pi\cdot\frac\delta{4\pi}=\frac\delta2.\end{equation} The statement of the proposition follows by adding this bound to the inequaltiy $\eqref{IntervalVarPropEq1}<\frac\delta2$.
\end{proof}

\section{Reduction to finitely many scales}
In order to show  for the map $T$ on $\bT^2$, it suffices to consider test functions $e_{(\zeta_1,\zeta_2)}(x,y)=e(\zeta_1x+\zeta_2y)$, $\zeta_1,\zeta_2\in\bZ$, since they span a dense subspace in $C^0(\bT^2)$. By taking complex conjugate, we can also assume $\zeta_2\geq 0$. From now on, fix such a function $f(x,y)=e_{(\zeta_1,\zeta_2)}(x,y)$.

The aim is to show that for any given $\delta>0$, for suffciently large $N$,
\begin{equation}\label{CorrEq}\left|\Exp_{n=0}^{N-1}\mu(n)(f\circ T^n)(x,y)\right|\end{equation} is bounded by $\delta$ for all $x$ and $y$.

Define $I_1^{k_+}$ according to Definition \ref{OstrowIntervalDef} and notice that it coincides with $\{0,\cdots, q_{k_++1}-1\}$. We decompose the correlation average between $f$ and $\mu$ into short averages along intervals of length $q_{k_++1}-1$. In fact,
\begin{equation}\label{DecomposedCorrEq1}\left|\Exp_{n=0}^{N-1}\mu(n)(f\circ T^n)(x,y)-\Exp_{L=0}^{N-1}\Exp_{n=1}^{q_{k_++1}-1}\mu(L+n)(f\circ T^{L+n})(x,y)\right|\ll\frac{q_{k_++1}}N,\end{equation} since the two averages only differ at worst by an $O(\frac{q_{k_++1}}N)$-portion of elements.

Together with \eqref{DecomposedCorrEq1},  Lemma \ref{IndepLem2} shows
\begin{equation}\label{DecomposedCorrEq2}\begin{aligned}\eqref{CorrEq}\leq&\Exp_{L=0}^{N-1}\left|\Exp_{r\in I_1^{k_--1}}\Exp_{s\in I_{k_-}^{k_+}}\mu(L+r+s)(f\circ T^{L+r+s})(x,y)\right|\\
&\hskip2cm+O\Big(\frac{q_{k_++1}}N+\frac1{a_{k_-}}\Big).\end{aligned}\end{equation}

One can write $(f\circ T^{L+r+s})(x,y)$ as
\begin{equation}\label{fExpressionEq}\begin{aligned}
&e_{(\zeta_1,\zeta_2)}\Big(x+(L+r+s)\alpha,\\
&\hskip2cm y+H_L(x)+H_r(x+L\alpha)+H_s(x+(L+r)\alpha)\Big)\\
=&e\Big(\zeta_1(r+s)\alpha+\zeta_2H_r(x+L\alpha)+\zeta_2H_n(x+(L+r)\alpha)\Big)\cdot e(\beta_{x,L}),
\end{aligned}\end{equation}
where $\beta_{x,L}$ is independent of $n$ and $r$. Therefore
\begin{equation}\label{SubCorrEq1}\begin{aligned}&\left|\Exp_{r\in I_1^{k_--1}}\Exp_{s\in I_{k_-}^{k_+}}\mu(L+r+s)(f\circ T^{L+r+s})(x,y)\right|\\
=&\Bigg|\Exp_{r\in I_1^{k_--1}}\Exp_{s\in I_{k_-}^{k_+}}\mu(L+r+s)\\
&\hskip1cm\cdot e\Big(\zeta_1(r+s)\alpha+\zeta_2H_r(x+L\alpha)
+\zeta_2H_n(x+(L+r)\alpha\Big)\Bigg|
\end{aligned}\end{equation}

\begin{lemma}\label{AlmostConstLem}For all $\delta>0$, if $k_-$ is large enough, then for all $k_+\geq k_-$ and all $x\in \bT^1$, $L\in\bN$, $$
\begin{aligned}&\left|\Exp_{r\in I_1^{k_--1}}\Exp_{s\in I_{k_-}^{k_+}}\mu(L+r+s)(f\circ T^{L+r+s})(x,y)\right|\\
\leq &\left|\Exp_{r\in I_1^{k_--1}}\Exp_{s\in I_{k_-}^{k_+}}\mu(L+r+s)e\Big(\zeta_1(r+s)\alpha+\zeta_2H_r(x+L\alpha)\Big)\right|+\frac\delta{8}.\end{aligned}$$
\end{lemma}
\begin{proof} By Proposition \ref{IntervalVarProp}, when $k_3\geq k_2(\frac\delta{16\zeta_2 +16})$, there is a constant $z_1$ that is independent of $x$, $L$, $r$ and $s$, such that $|z_1|\leq 1$, and with probability $1-\frac\delta{16}$ or higher for a randomly chosen $s\in I_{k_-}^{k_+}$, \begin{equation}\label{AlmostConstLemEq1}\left|e\Big(\ksum \phi_k(s_k)\Big)-z_1\right|<\frac\delta{16\zeta_2 +16}.\end{equation}

We write $x+(L+r)\alpha\in\bT^1$ as $\sum_{k=1}^\infty(\tx_k+L_k+r_k)\theta_k$. Remark that by construction, $r_k=0$ for $k\geq k_-$. Then, again with probability at least $1-\frac\delta{16}$, \begin{equation}\label{AlmostConstLemEq2}\left|e\Big(\ksum \phi_k\big(s_k+(\tx_k+L_k)\big)-\phi_k(\tx_k+L_k)\Big)-z_2\right|<\frac\delta{16\zeta_2 +16},\end{equation} where $z_2=z_1e(-\phi_k(\tx_k+L_k))$ depends on $x$ and $L$ but not on $r$ and $n$.
Here we used the fact that $\phi_k$ is a function defined on $\bZ/a_k\bZ$.

By Proposition \ref{RandomSumProp}, \eqref{AlmostConstLemEq2} becomes
\begin{equation}\label{AlmostConstLemEq3}
\begin{aligned}\left|e\big(H_s(x+(L+r)\alpha)\big)-z_2\right|<&\frac\delta{16\zeta_2 +16}+\frac\delta{16\zeta_2 +16}\\
=&\frac\delta{16\zeta_2 +16},\end{aligned}\end{equation} if we choose $k_3\geq k_0(\frac1{2\pi}\cdot\frac\delta{16\zeta_2 +16})$. Note that Proposition \ref{RandomSumProp} applies because $\{n_k\}$ is an Ostrowski numeration, and $|\tx_k+L_k+r_k|\leq 3|a_k|$ for all $k$.

Remark a simple fact: if $|w_1|, |w_2|\leq 1$, then $$|w_1^{\zeta_2}-w_2^{\zeta_2}|=|w_1-w_2|\cdot\big|\sum_{j=0}^{\zeta_2}w_1^jw_2^{\zeta_2-j}\big|\leq (\zeta_2+1)\cdot|w_1-w_2|.$$ So
\begin{equation}\label{AlmostConstLemEq16}
\left|e\big(\zeta_2H_s(x+(L+r)\alpha)\big)-z_2^{\zeta_2}\right|<(\zeta_2+1)\cdot\frac\delta{16\zeta_2 +16}<\frac\delta{16}.\end{equation}
Remember that this does not necessarily hold always, but with a probability higher than $1-\frac\delta{16\zeta_2 +16}$ with respect to $s\in I_{k_-}^{k_+}$.

Plugging this estimate into \eqref{SubCorrEq1}, we see that, thanks to the fact that $z_2^{\zeta_2}$ doesn't depend on $r$ or $n$,
\begin{equation}\label{AlmostConstLemEq4}
\begin{aligned}&\Bigg|\Exp_{r\in I_1^{k_--1}}\Exp_{s\in I_{k_-}^{k_+}}\mu(L+r+s)(f\circ T^{L+r+s})(x,y)\\
&\hskip5mm-z_2^{\zeta_2}\Exp_{r\in I_1^{k_--1}}\Exp_{s\in I_{k_-}^{k_+}}\mu(L+r+s)e\Big(\zeta_1(r+s)\alpha+\zeta_2H_r(x+L\alpha)\Big)\Bigg|\\
<&\frac\delta{16\zeta_2 +16}+\frac\delta{16}<\frac\delta{8}.\end{aligned}\end{equation}

We obtain the lemma because $|z_2|\leq 1$ and $\zeta_2\geq 0$.
\end{proof}

For each $n\geq 0$, truncate its Ostrowski numeration before the $k_-$-th digit, defining a function $r=r(n)=\sum_{k=1}^{k_--1}n_kq_k\in I_1^{k_--1}$.

\begin{corollary}\label{AlmostConstCor}For all $\delta>0$, there exists $k_3=k_3(\alpha,\delta,\zeta_1,\zeta_2)\in\bN$, such that:

For all $k_+\geq k_-\geq k_3$, and all $x\in\bT^1$,
$$\begin{aligned}&\Exp_{n=0}^{N-1}\mu(n)(f\circ T^n)(x,y)\\
\leq&\Exp_{L=0}^{N-1}\left|\Exp_{n\in I_1^{k_+}}\mu(L+n)e\big(\zeta_1n\alpha+\zeta_2H_{r(n)}(x+L\alpha)\big)\right|+\frac{\delta}{8}+O\Big(\frac{q_{k_++1}}N+\frac1{a_{k_-}}\Big)\end{aligned}$$
\end{corollary}

\begin{proof}This follows from the inequality \eqref{DecomposedCorrEq2}, Lemma \ref{AlmostConstLem}, and one more application of Lemma \ref{IndepLem2}.\end{proof}

\section{Proof of the main theorem}

We can think of $r(n)$ as a quasi-periodic residue of $n$. Indeed, from the construction of the Ostrowski numeration, one can easily see the following fact.

\begin{lemma}\label{ResidueLem} For each $r\in I_1^{k_--1}=\{0,\cdots,q_{k_-}-1\}$, there is an arc\footnote{Here an arc can be either open or closed at each of its endpoints.} $D_r\subset\bT^1$ such  that:\begin{enumerate}                                                                                                       \item $\bT^1$ decomposes as a disjoint union $\bigsqcup_{r=0}^{q_{k_-}-1}D_r$;
\item For $n\geq 0$, $r(n)=r$ if and only if $n\alpha\in D_r$ modulo $\bZ$.\end{enumerate}
\end{lemma}

The lemma allows to interprete $e\big(\zeta_1n\alpha+\zeta_2H_{r(n)}(x+L\alpha)\big)$ as a function of $n\alpha$.

\begin{proposition}\label{ExpDecompProp} For all $\delta>0$ and $f(x,y)=e_{(\zeta_1,\zeta_2)}(x,y)$ , there exists $k_3=k_3(\alpha,\delta,\zeta_1,\zeta_2)\in\bN$, such that:

For all $k_-\geq k_3$ , there exist $A=A(\alpha,  \delta, k_-)$, $B=B(\alpha,  \delta, k_-)$and $k_4=k_4(\alpha,\delta,k_-)$, such that for all $(x,y)\in\bT^2$ and $k_+\geq k_4$,
$$\begin{aligned}&\left|\Exp_{n=0}^{N-1}\mu(n)(f\circ T^n)(x,y)\right|\\
\leq& B\cdot \max_{\xi=-A}^A \Exp_{L=0}^{N-1}\left|\Exp_{n=0}^{q_{k_+}-1}\mu(L+n)e\big((\zeta_1+\xi)n\alpha\big)\right| +\frac{\delta}{4}+O\Big(\frac{q_{k_++1}}N+\frac1{a_{k_-}}\Big).\end{aligned}$$ \end{proposition}
\begin{proof} Suppose $k_-$ is fixed, for every $r=0, \cdots,q_{k_-}-1$, choose a continuous function $\psi_r:\bT^1\mapsto[0,1]$ approximating $\bfone_{D_r}$ in the sense that $\psi_r=\bfone_{D_r}$ except in two short arcs $U_r^-$ and $U_r^+$ respectively around both ends of $D_r$ with total length $|U_r^-|+|U_r^+|<\frac\delta{32q_{k_-}}$.

Furthermore, because trigonometric polynomials are dense in $C^0(\bT^1)$, one can find $A_r\in\bN$ such that $\psi_r(w)$ is approximated by a trigonometric polynomial $\sum_{\xi=-A_r}^{A_r}\theta_{r,\xi}e(\xi w)$ up to error $\frac\delta{16q_{k_-}}$ in $C^0$ norm. By taking maximum over all $r$, we can make $A_r=A$ a constant that is determined by $\alpha$, $\delta$ and $k_-$.

When $L$ and $x$ are given, the function
\begin{equation}\label{ExpDecompPropEq1}\sum_{r=0}^{q_{k_-}-1}e(\zeta_2H_r(x+L\alpha))\bfone_{D_r}(w)\end{equation}
is approximated by
\begin{equation}\label{ExpDecompPropEq2}\sum_{\xi=-A}^Ac_{\xi,L,x}e(\xi w):=\sum_{r=0}^{q_{k_-}-1}e(\zeta_2H_r(x+L\alpha))\sum_{\xi=-A}^A\theta_{r,\xi}e(\xi  w)\end{equation}
up to an error of $q_{k_-}\cdot \frac\delta{16q_{k_-}}=\frac\delta{16}$, unless $w$ lies in the exceptional set $U=\bigcup_{r=0}^{q_{k_-}-1}(U_r^-\cup U_r^+)$.

Let $B_0=\max_{\xi=-A}^A\sum_{r=0}^{q_{k_-}-1}|\theta_{r,\xi}|$, then $B_0$ is determined by $\alpha$, $\delta$ and $k_-$, and satisfies
\begin{equation}\label{ExpDecompPropEq3}\max_{\xi=-A}^A|c_{\xi,L,x}|\leq B_0, \forall L\geq 0,\forall x\in\bT^1.\end{equation}

By Lemma \ref{ResidueLem}, $e(\zeta_2H_{r(n)}(x+L\alpha))=\sum_{r=0}^{q_{k_-}-1}e(\zeta_2H_r(x+L\alpha))\bfone_{D_r}(n\alpha)$. Therefore,
$$\begin{aligned}&\left|\mu(L+n)e\big(\zeta_1n\alpha+\zeta_2H_{r(n)}(x+L\alpha)\big)
-\sum_{\xi=-A}^Ac_{\xi,L,x}\mu(L+n)e\big((\zeta_1+\xi)n\alpha\big)\right|\\
\leq&\frac\delta{16}\end{aligned}$$
unless $n\alpha\in U$.

Because $|U|=\frac\delta{32}$ and the sequence $\{n\alpha\}$ is equidistributed in $\bT^1$, for some $Q_0$ that depends on $U$, $\Prob_{0\leq n\leq Q-1}(n\alpha\in U)\leq\frac\delta{16}$ for all $Q\geq Q_0$. Choose the smallest $k_4$ such that $q_{k_4}>Q_0$, then $k_4$ is determined by $\alpha$, $\delta$ and the choice of $U$, which in turn depends on $k_-$. Therefore, for all $L$, $x$ and $k_+\geq k_4$,

\begin{equation}\label{ExpDecompPropEq4}\begin{aligned}&\left|\Exp_{n=0}^{q_{k_+}-1}\mu(L+n)e\big(\zeta_1n\alpha+\zeta_2H_{r(n)}(x+L\alpha)\big)\right.\\
&\hskip1cm\left.-\Exp_{n=0}^{q_{k_+}-1}\sum_{\xi=-A}^Ac_{\xi,L,x}\mu(L+n)e\big((\zeta_1+\xi)n\alpha\big)\right|\\
\leq&\frac\delta{16}+\frac\delta{16}=\frac\delta{8}.\end{aligned}\end{equation}

Plugging \eqref{ExpDecompPropEq4} into the inequality in Corollary \ref{AlmostConstCor} yields that

$$\begin{aligned}&\Exp_{n=0}^{N-1}\mu(n)(f\circ T^n)(x,y)\\
\leq&\Exp_{L=0}^{N-1}\left|\Exp_{n=0}^{q_{k_+}-1}\sum_{\xi=-A}^Ac_{\xi,L,x}\mu(L+n)e\big((\zeta_1+\xi)n\alpha\big)\right|+\frac{\delta}{4}+O\Big(\frac{q_{k_++1}}N+\frac1{a_{k_-}}\Big)\\
\leq&\Exp_{L=0}^{N-1}\sum_{\xi=-A}^A|c_{\xi,L,x}|\cdot\left|\Exp_{n=0}^{q_{k_+}-1}\mu(L+n)e\big((\zeta_1+\xi)n\alpha\big)\right|+\frac{\delta}{4}+O\Big(\frac{q_{k_++1}}N+\frac1{a_{k_-}}\Big)\\
\leq& B_0\cdot \sum_{\xi=-A}^A \Exp_{L=0}^{N-1}\left|\Exp_{n=0}^{q_{k_+}-1}\mu(L+n)e\big((\zeta_1+\xi)n\alpha\big)\right|+\frac{\delta}{4}+O\Big(\frac{q_{k_++1}}N+\frac1{a_{k_-}}\Big).
\end{aligned}$$
The proposition follows by setting $B=B_0\cdot (2A+1)$.
\end{proof}

The proof of Theorem \ref{MainThm}  relies on the recent theorem of Matom\"aki, Radziwi\l\l, and Tao \cite{MRT15} on averages of non-pretentious multiplicative functions in short intervals, or more precisely the following corollary of it:
\begin{proposition}\label{MRTProp}\cite{MRT15}  For some constant $c$, for all $N\geq R\geq 10$ and $\beta\in\bR$,
$$\Exp_{L=0}^{N-1}\left|\Exp_{n=1}^R\mu(L+n)e(\beta n)\right|\ll(\log N)^{-c}+\frac{\log\log R}{\log R}.$$ The implied constant does not depend on $\beta$.\end{proposition}
This is \cite{MRT15}*{Theorem 1.7}, applied to $\mu$ in light of the paragraph before the theorem in that paper, which asserts that $\mu$ is sufficiently non-pretentious.

\begin{proof}[Proof of Theorem \ref{MainThm}]

We now give the proof of our main result. Recall that it suffices to work under Hypothesis \ref{ReducedHypo} and assume that $f(x,y)=e_{(\zeta_1,\zeta_2)}(x,y)$. In this case Proposition \ref{ExpDecompProp} applies.

Using Proposition \ref{MRTProp}, the inequality in Proposition \ref{ExpDecompProp} becomes
\begin{equation}\label{MainPfEq1}\begin{aligned}&\left|\Exp_{n=0}^{N-1}\mu(n)(f\circ T^n)(x,y)\right|\\
\leq&O\Big(B (\log N)^{-c}+B\cdot\frac{\log\log q_{k_++1}}{\log q_{k_++1}}\Big) +\frac\delta{4}+O\Big(\frac{q_{k_++1}}N+\frac1{a_{k_-}}\Big).\end{aligned}\end{equation}

Given $\alpha$, $\delta$, $\zeta_1$, $\zeta_2$, we can first choose a large $k_-\geq k_3(\alpha,\delta,\zeta_1,\zeta_2)$ which satisfies the extra condition that $q_{k_-+1}> q_{k_-}^{\frac\tau2}$. Recall that there are infinitely many such $k_-$ under Hypothesis \ref{ReducedHypo}. Then by \eqref{akGrowthEq}, one may fix a sufficiently large $k_-$ for which the $O(\frac 1{a_{k_-}})$ term is bounded by $\frac\delta 8$. Note $B$ is determined once $k_-$ is fixed. One can then pick a large $k_+\geq k_4(\alpha,\delta,k_-)$ such that the $O \big(B\cdot \frac{\log\log q_{k_++1}}{\log q_{k_++1}})$ part is bounded by $\frac\delta8$. Finally, in order to make the entire sum in \eqref{MainPfEq1} less than $\delta$, it suffices to make $N$ much larger compared with $q_{k_+}$ and $B$ to make the two other terms arbitrarily small. Since $\delta$ is arbitrary, the proof is completed.
\end{proof}


\begin{bibdiv}
\begin{biblist}

\bib{B01}{article}{
   author={Berth{\'e}, Val{\'e}rie},
   title={Autour du syst\`eme de num\'eration d'Ostrowski},
   language={French, with French summary},
   journal={Bull. Belg. Math. Soc. Simon Stevin},
   volume={8},
   date={2001},
   number={2},
   pages={209--239},
}

\bib{B13a}{article}{
   author={Bourgain, J.},
   title={On the correlation of the Moebius function with rank-one systems},
   journal={J. Anal. Math.},
   volume={120},
   date={2013},
   pages={105--130},
}

\bib{B13b}{article}{
   author={Bourgain, J.},
   title={M\"obius-Walsh correlation bounds and an estimate of Mauduit and
   Rivat},
   journal={J. Anal. Math.},
   volume={119},
   date={2013},
   pages={147--163},
}

\bib{BSZ13}{article}{
   author={Bourgain, J.},
   author={Sarnak, P.},
   author={Ziegler, T.},
   title={Disjointness of M\"obius from horocycle flows},
   conference={
      title={From Fourier analysis and number theory to Radon transforms and
      geometry},
   },
   book={
      series={Dev. Math.},
      volume={28},
      publisher={Springer, New York},
   },
   date={2013},
   pages={67--83},
}

\bib{D37}{article}{
   author={Davenport, H.},
   title={On some infinite series involving arithmetical functions II},
   journal={Quat. J. Math.},
   volume={8},
   date={1937},
   pages={313--320},
}

\bib{ELD14}{article}{
   author={El Abdalaoui, El Houcein},
   author={Lema{\'n}czyk, Mariusz},
   author={de la Rue, Thierry},
   title={On spectral disjointness of powers for rank-one transformations
   and M\"obius orthogonality},
   journal={J. Funct. Anal.},
   volume={266},
   date={2014},
   number={1},
   pages={284--317},
}

\bib{FKLM15}{article}{
   author={Ferenzi, S},
   author={Ku\l{}aga-Przymus, Joanna},
   author={Lema\'nczyk, Mariusz},
   author={Mauduit, C},
   title={Substitutions and Möbius disjointness},
   journal={preprint},
   date={2015},
}

\bib{F61}{article}{
   author={Furstenberg, H.},
   title={Strict ergodicity and transformation of the torus},
   journal={Amer. J. Math.},
   volume={83},
   date={1961},
   pages={573--601},
}

\bib{G12}{article}{
   author={Green, Ben},
   title={On (not) computing the M\"obius function using bounded depth
   circuits},
   journal={Combin. Probab. Comput.},
   volume={21},
   date={2012},
   number={6},
   pages={942--951},
}

\bib{GT12}{article}{
   author={Green, Ben},
   author={Tao, Terence},
   title={The M\"obius function is strongly orthogonal to nilsequences},
   journal={Ann. of Math. (2)},
   volume={175},
   date={2012},
   number={2},
   pages={541--566},
}

\bib{KL15}{article}{
   author={Ku{\l}aga-Przymus, J.},
   author={Lema{\'n}czyk, M.},
   title={The M\"obius function and continuous extensions of rotations},
   journal={Monatsh. Math.},
   volume={178},
   date={2015},
   number={4},
   pages={553--582},
}


\bib{K97}{book}{
   author={Khinchin, A. Ya.},
   title={Continued fractions},
   edition={Translated from the third (1961) Russian edition},
   note={With a preface by B. V. Gnedenko;
   Reprint of the 1964 translation},
   publisher={Dover Publications, Inc., Mineola, NY},
   date={1997},
   pages={xii+95},
}

\bib{LS15}{article}{
   author={Liu, Jianya},
   author={Sarnak, Peter},
   title={The M\"obius function and distal flows},
   journal={Duke Math. J.},
   volume={164},
   date={2015},
   number={7},
   pages={1353--1399},
   issn={0012-7094},
}

\bib{MMR14}{article}{
   author={Martin, Bruno},
   author={Mauduit, Christian},
   author={Rivat, Jo{\"e}l},
   title={Th\'eor\'eme des nombres premiers pour les fonctions digitales},
   language={French},
   journal={Acta Arith.},
   volume={165},
   date={2014},
   number={1},
   pages={11--45},
}

\bib{MRT15}{article}{
   author={Matom\"aki, Kaisa},
   author={Radziwi\l\l, Maksym},
   author={Tao, Terence},
   title={An averaged form of Chowla's conjecture},
   journal={preprint},
   date={2015},
}

\bib{MR10}{article}{
   author={Mauduit, Christian},
   author={Rivat, Jo{\"e}l},
   title={Sur un probl\`eme de Gelfond: la somme des chiffres des nombres
   premiers},
   language={French, with English and French summaries},
   journal={Ann. of Math. (2)},
   volume={171},
   date={2010},
   number={3},
   pages={1591--1646},
}

\bib{O22}{article}{
   author={Ostrowski, Alexander},
   title={Bemerkungen zur Theorie der Diophantischen Approximationen},
   language={German},
   journal={Abh. Math. Sem. Univ. Hamburg},
   volume={1},
   date={1922},
   number={1},
   pages={77--98},
}

\bib{P15}{article}{
   author={Peckner, Ryan},
   title={M\"obius disjointness for homogeneous dynamics},
   journal={preprint},
   date={2015},
}

\bib{S09}{article}{
   author={Sarnak, Peter},
   title={Three lectures on the M\"obius function, randomness and dynamics},
   journal={lecture notes, IAS},
   date={2009},
}

\bib{V77}{article}{
   author={Vaughan, Robert-C.},
   title={Sommes trigonom\'etriques sur les nombres premiers},
   language={French, with English summary},
   journal={C. R. Acad. Sci. Paris S\'er. A-B},
   volume={285},
   date={1977},
   number={16},
}

\bib{W15}{article}{
   author={Wang, Zhiren},
   title={M\"obius disjointness for analytic skew products},
   journal={preprint},
   date={2015},
}

\end{biblist}
\end{bibdiv}
\end{document}